\documentclass[reqno]{amsart}

\usepackage{amssymb}
\usepackage{amsmath}
\usepackage{verbatim}
\newtheorem{theorem}{Theorem}
\newtheorem{lemma}{Lemma}
\newtheorem{proposition}{Proposition}

\newtheorem{definition}{Definition}

\newtheorem{corollary}{Corollary}
\newtheorem{remark}{Remark}
\newtheorem{assumption}{Assumption}

\begin{document}
\title{On a time and space discretized approximation of the
Boltzmann equation in the whole space}
\author{C. P. Gr\"unfeld$^\ast$$^\dag$}\thanks{$^\ast$
Institute of Space Science, PO Box MG-23, Magurele - Ilfov, Romania
}
\author{D. Marinescu$^\dag$}\thanks{$^\dag$
Institute
of Mathematical Statistics and Applied Mathematics of the
Romanian  Academy,
PO Box 1-24, Bucharest, Romania
}
\thanks{
{\sc Email  addresses:} {\bf grunfeld@spacescience.ro} (C. P.  Gr\"unfeld)
{\bf dorin.marinescu@ima.ro} (D. Marinescu)} 

\begin{abstract}
  In this paper, convergence results on the solutions of a time and space
  discrete model approximation of the Boltzmann equation for a gas of
  Maxwellian particles in a bounded domain, obtained by Babovsky and Illner
  [1989], are extended to  approximate the solutions of the Boltzmann
  equation in the whole physical space. This is done for a class of particle
  interactions including Maxwell and soft cut-off potentials in the sense of
  Grad.

  The main result shows that the solutions of the discrete model converge
  in $\mathbb{L}^1$ to the solutions of the Boltzmann equation, when the
 discretization parameters  go simultaneously to zero. The convergence is
  uniform with respect to the  discretization parameters.

  In addition, a sufficient condition for the implementation
  of the main result is provided.
\end{abstract}
\maketitle
\pagestyle{myheadings}
\markboth{On a time and space discretized approximation of
  the Boltzmann equation in the whole space}{C. P. Gr\"unfeld and D. Marinescu}

\section{Introduction \label{s1}}

In a known paper \cite{BBb}, Babovsky and Illner provided a validation
(convergence) proof of Nanbu's simulation method \cite{NB} for the spatially
inhomogeneous (full) Boltzmann equation \cite{LB,Cer} describing a rarefied
gas of Maxwellian particles confined to a bounded spatial domain (with
specularly reflecting boundary conditions). More specifically, the main
result (Theorem 7.1) of \cite{BBb} demonstrated that the discrete measures
provided by Nanbu's simulation method are almost surely weakly convergent to
absolutely continuous measures with densities given by solutions of the
Boltzmann equation.

In essence, the analysis behind the main theorem of
\cite{BBb} represented a space-dependent generalization of
the convergence proof of Nanbu's simulation algorithm for
the space-homogeneous Boltzmann equation, provided in an
earlier work by Babovsky \cite{Bab89}. Briefly, in
\cite{BBb}, the space-homogeneous simulation algorithm of
\cite{Bab89} was applied to a suitable time and space
discrete Boltzmann model (Eq. (5.14) in \cite{BBb}).  The
latter was derived from the Boltzmann equation by means of
time discretization, splitting (separation of free flow and
collisional interactions), cell-partitioning of the physical
domain of the gas, and space-averaging (homogenization) over
cells. The discretization was parameterized by a time-step
and an upper bound for the maximum of all cell
diameters. The analysis was completed by combining
convergence properties of the discrete Boltzmann model with
those of the space-homogeneous simulation algorithm of
\cite{Bab89}.  To this end, Babovsky and Illner established
the key result (Corollary 5.1 in \cite{BBb}) that the
solutions of the discrete Boltzmann model converge in
discrepancy\footnote {Let $\mu$ and $\nu$ be two (Borel)
  probability measures on the same measure space
  $\mathfrak{B}\subseteq \mathbb{R}^n $.  Consider
  $\mathfrak{B}$ with the usual semi-order $\leq$ of
  $\mathbb{R}^n $. Following \cite{BBb}, the discrepancy
  between $\mu$ and $\nu$ is defined as
  $D(\mu,\nu):=sup_{z\in \mathbb{R}^n}|\mu\{y\in
  \mathfrak{B}: y \leq z \}-\nu\{y\in \mathfrak{B}: y \leq z \}|
  $.} to the solutions of the Boltzmann equation for the gas
in a bounded spatial volume, uniformly with respect to the
parameters of the discretization, when these parameters
converge simultaneously to zero.

A notable thing about the proof of the convergence in
discrepancy of the solutions of the discretized Boltzmann
model is that, as it appears in \cite{BBb}, is responsible
for the limitation of the analysis of \cite{BBb} to the case
of the Boltzmann gas in a bounded spatial domain. Indeed,
the boundedness of the spatial domain was actually assumed
in \cite{BBb} in order to prove the above key convergence
result for the discrete Boltzmann equation (see \cite{BBb},
p. 59). An alternative proof, without the boundedness
assumption, might allow the conclusions of \cite{BBb} to be
extended to other important examples, e.g., a gas expanding
in the whole physical space.

In this paper, the results on the convergence in discrepancy established by
Babovsky and Illner for the solutions of their discrete Boltzmann model of
\cite{BBb} are extended  to the setting of the Boltzmann equation in the
entire physical space. More specifically, in such a setting, we show that the
solutions of the discrete Boltzmann model converge in $\mathbb{L}^1$ to the
solutions of the Boltzmann equation in the whole physical space, uniformly
with respect to the parameters of the discretization, when these parameters
converge simultaneously to zero. We also show that the solutions of the
discrete approximation satisfy the conservation laws for mass, momentum and
energy.

Here, it should be recalled that the results of \cite{BBb}
concern the Boltzmann equation for Maxwellian particles. The
limitation to Maxwellian interactions does not come from
the proof of the analytical convergence of the discretized
Boltzmann model, but is imposed by the implementation of the
simulation algorithm of \cite{Bab89} for the validation of
Nanbu's scheme (see \cite{BBb}, p. 48).  However, besides
its usefulness in the validation of the Nanbu's scheme, the
discrete Boltzmann model of \cite{BBb} might be applied to
obtain new (not necessarily probabilistic) rigorous
algorithms for the Boltzmann equation.  Thus, understanding
its convergence properties in more general situations than
in \cite{BBb} may be of interest. In this respect, as an
additional contribution, our main result concerns the
Boltzmann equation with Maxwell and soft cut-off collision
kernels in the sense of Grad \cite{Grad}.

Compared to \cite{BBb}, our analysis must face additional difficulties, since
one has to estimate, uniformly, in some sense, how high speed gas particles
situated at large distances contribute to the gas evolution, approximated as
in \cite{BBb}, by an alternation of molecular transport and collision steps.
In this respect, a technical point is reconsidering the important property
established by Babovsky and Illner (Theorem~5.1 in \cite{BBb}) that, under
suitable conditions, if the Boltzmann equation is approximated by the
discrete Boltzmann model, then the family of errors introduced by the
approximation is bounded in some $\mathbb{L}^{\infty}$ - (velocity)
Maxwellian weighted space, uniformly with respect to the parameters of the
discretization. This property was demonstrated in \cite{BBb}, in the setting
of the Boltzmann equation in a bounded spatial domain, but remains actually
valid in a larger context, as is implicit from \cite{BBb}. Nevertheless, for
the sake of clarity and completeness, in the present work, we will prove a
precise statement appropriate to our framework (see
Proposition~\ref{est-del-tx-BI} in Subsection~\ref{s4.2}).

The rest of this paper is structured as follows. In Section ~\ref{s1b}, we
present the discrete Boltzmann model of \cite{BBb}, and formally introduce
the main result. However, a precise formulation (Theorem~\ref{conv-discr}) is
given in the second part of Section~\ref{s2}. This requires some preparation
in the first part of the same section. The second part of Section~\ref{s2}
also includes Theorem \ref{spatial-lip} which provides sufficient conditions
for the application of Theorem~\ref{conv-discr}. Section~\ref{s3} deals with
the proofs of the theorems stated in Section ~\ref{s2}. The proofs rely on
technical estimates  provided in Subsection~\ref{s4.1}. In particular,
standard $\mathbb{L}^{\infty}$ - type inequalities for the collision term are
adapted to our setting, supplemented with useful $\mathbb{L}^{1}$ -
estimates.  The central result of Subsection~\ref{s4.1} is Lemma~\ref{key},
needed later to measure, in some sense, the errors introduced when the
discrete Boltzmann model approximates  the Boltzmann equation. The results of
Subsection~\ref{s4.1} are then used in Subsection~\ref{s4.2} to prove
Proposition~\ref{est-del-tx-BI}, ultimately leading to the proof of
Theorem~\ref{conv-discr}. Subsection~\ref{s4.3} contains the proof of
Theorem~\ref{spatial-lip}. Finally, Section~\ref{s5} presents a simple
application to a Boltzmann model for a rarefied gas expanding in the whole
space, and closes with a few concluding remarks and possible future
directions.

\section{Discretized Boltzmann model for the Boltzmann
  equation \label{s1b}} In this section we recall some very
basic facts about the Boltzmann equation, and briefly
present its time and space discrete approximation of
\cite{BBb}, adapted to our setting. Finally, we formally
introduce our results.

The Cauchy problem for the Boltzmann equation for a simple
gas (monatomic gas of identical particles with elastic
binary collisions), evolving in the whole physical space
reads (in non-dimensional units) as \cite{LB,Cer}
\begin{equation}
\begin{array}{l} \displaystyle
\frac{\partial f}{\partial
    t}+\mathbf{v}\cdot \nabla _{\mathbf{x} }f=J(f) \quad
  {\rm in}
  \quad  (0,T)\times \mathbb{R}^{3} \times \mathbb{R}^{3},
  \\ \\    \displaystyle
f(0,\mathbf{x},\mathbf{v})=f_{0}(\mathbf{x},\mathbf{v})
\quad {\rm on} \quad  \mathbb{R}^3 \times \mathbb{R}^3,
\end{array}
\label{ec-bg}
\end{equation}
where the unknown $f=f(t,\mathbf{x},\mathbf{v})\geq 0$
represents the distribution density of the gas particles at
time $0 \leq t <T\leq \infty $, with position
$\mathbf{x}=(x_{1},x_{2},x_{3}) \in \mathbb{R}^{3}$, and
velocity $\mathbf{v}=(v_{1},v_{2},v_{3}) \in
\mathbb{R}^{3}$. The right hand side of the above equation
is the nonlinear Boltzmann collision term
$J(f):=J_{B}(f,f)$, where
\begin{equation}
  J_{B}(g,h)(t,\mathbf{x},\mathbf{v})=\int\limits_{\mathbb{R}^{3}\times
    \mathbb{S}^{2}}b(|\mathbf{v }-\mathbf{v}_{\ast
  }|,\boldsymbol{\omega })(g^{\prime }h_{\ast }^{\prime
  }-gh_{\ast })d\mathbf{v}_{\ast }d\boldsymbol{\omega
  }. \label{jgh}
\end{equation}  acts on $g$ and $h$ only through the
velocity dependence. Here, we are using the standard shorthand notations
$g:=g(t, \mathbf{x},\mathbf{v})$, $h_{\ast }:=h(t,
\mathbf{x},\mathbf{v}_{\ast})$, $g^{\prime }:=g(t,
\mathbf{x},\mathbf{v}^{\prime })$, $h_{\ast }^{\prime }:=h(t,
\mathbf{x},\mathbf{v}_{\ast }^{\prime })$. Moreover,
\begin{equation} \mathbf{v}^{\prime
  }=\mathbf{v}-((\mathbf{v}-\mathbf{v}_{\ast
  })\cdot\boldsymbol{\omega })\boldsymbol{\omega }, \quad
\mathbf{v}_{\ast }^{\prime }=\mathbf{v} _{\ast
  }+((\mathbf{v}-\mathbf{v}_{\ast })\cdot\boldsymbol{\omega
  })\boldsymbol{ \omega },
  \label{param}
  \end{equation}
  are the post-collisional velocities expressed, in terms of
  the pre-collisional velocities ($\mathbf{v}$,
  $\mathbf{v}_{\ast })$ and the collision parameter
  $\boldsymbol{\omega}$ over the unit sphere
  $\mathbb{S}^{2}$, as (parameterized) solutions of the laws
  of momentum and energy conservation in binary elastic
  collisions
\begin{equation} \mathbf{v}+\mathbf{v}_{\ast
  }=\mathbf{v}^{\prime }+\mathbf{v}_{\ast }^{\prime },\quad
  \mathbf{v}^{2}+\mathbf{v}_{\ast }^{2}=\mathbf{v}^{\prime
    2}+ \mathbf{v}_{\ast }^{\prime 2}.
    \label{econs}
\end{equation}
Furthermore, in (\ref{jgh}), the collision kernel
$b(|\mathbf{v}-\mathbf{v}_{\ast}|, \boldsymbol{\omega})$ is a given
nonnegative function that depends only on the modulus of the relative
velocity $|\mathbf{v}-\mathbf{v}_{\ast}|$ and the scalar product
$({\mathbf{v}-\mathbf{v}_{\ast}})\cdot \boldsymbol{\omega}$.

A useful decomposition of (\ref{jgh}) is $J_{B} (g,h)=P_{B}(g,h)-S_{B}(g,h)$
where
\begin{equation}
\begin{array}{l} \displaystyle
P_{B}(g,h)(\mathbf{x},\mathbf{v}):=\int\limits_{\mathbb{R}^{3}\times
    \mathbb{S}^{2}}b(|\mathbf{v}- \mathbf{v}_{\ast
  }|,\boldsymbol{\omega })g^{\prime }h_{\ast }^{\prime }d
  \mathbf{v}_{\ast }d\boldsymbol{\omega },
\\ \\ \displaystyle
 S_{B}(g,h)(\mathbf{x},\mathbf{v}):=\int\limits_{\mathbb{R}^{3}\times
    \mathbb{S}^{2}}b(|\mathbf{v}- \mathbf{v}_{\ast
  }|,\boldsymbol{\omega })gh_{\ast }d\mathbf{v}_{\ast }d
  \boldsymbol{\omega },
\end{array}
\label{ssb}
\end{equation}
provided that the integrals are well defined. We refer to \cite{Cer} for a
comprehensive presentation of the Boltzmann equation and its applications.
Here we just recall the property
\begin{equation} \int\limits_{\mathbb{R}^{3}}\varphi
  _{i}(\mathbf{v}
  )J(g)(\mathbf{x},\mathbf{v})d\mathbf{v}=0,\,i=0,1,2,3,4,
\label{intJ}
\end{equation}
valid for any $g$ for which the above integral exists, where
$\varphi _{0}(\mathbf{v}):=1,\varphi
_{i}(\mathbf{v}):=v_{i},\,i=1,2,3,\,\varphi_{4}(
\mathbf{v}):=\mathbf{v}^{2}$. Formally, this implies that
the solution of the Boltzmann equation satisfies the fluid
balance laws for mass, momentum and energy \cite{Cer}.

Let $0<\Delta t<T$ be a discretization time-step for the
time interval $[0,T]$. Suppose that $\mathbb{R}^{3}$ is
partitioned into a countable family $\Pi:=(\pi _{l})_{l\in
  \mathbb{N}}$ of distinct cells $\pi _{l}$ with finite
diameter $D(\pi _{l})$. Then the discrete model provided in \cite{BBb}
to approximate the Boltzmann equation in a bounded domain
(Eq. (5.14) in \cite{BBb}) can be reformulated in our setting (to approximate
(\ref{ec-bg})) as
\begin{equation} \left\{
\begin{array}{l}
  \tilde{f}^{0}(\mathbf{x},\mathbf{v})=
  f_{0}(\mathbf{x},\mathbf{v}), \\\\
  \tilde{f}^{j}(\mathbf{x},\mathbf{v})=U^{\Delta
    t}\tilde{f}^{j-1}(\mathbf{x}, \mathbf{v})  +\Delta
  t  J^{\pi}(U^{\Delta t}\tilde{f}^{j-1})(\mathbf{x},
  \mathbf{v}),\;\; j=1,\ldots,[[\frac{T}{\Delta t}]].
\end{array} \right. \label{discr3-BI}
\end{equation} Here, $U^{t}$ is the free streaming operator
\begin{equation}
  (U^{t}g)(\mathbf{x},\mathbf{v}):=g(\mathbf{x}-t\mathbf{v},\mathbf{v}),\quad
  t \in \mathbb{R}, \label{ut}
\end{equation} $J^{\pi}$ is the homogenized Boltzmann operator
\begin{equation} J^{\pi}(g):=J_{B}(g,{\pi}g), \label{jip0}
\end{equation}
defined by means of  the operator of homogenization of the space cells
\begin{equation} (\pi g)(\mathbf{x},\mathbf{v}):=\sum_{l\in
    \mathbb{N}}\chi_{\pi _{l}}(\mathbf{x})
  \frac{1}{| \pi _{l}|
  }\int\limits_{\pi_{l}}g(\mathbf{y},
  \mathbf{v})\,d\mathbf{y},
\label{pi-omo}
\end{equation}
where $\chi _{\pi_{l}}$ is the indicator function of the
cell $\pi _{l}$, $|\pi_{l}|$ denotes the volume of the cell.
Moreover, $[[\cdot]]$ is the integer part function.

Apart from minor changes of notation, the main difference
between (\ref{discr3-BI}) and the formulation of the
discretized Boltzmann model of \cite{BBb} is that the sum in
(\ref{discr3-BI}) is infinite, because, in our setting, we
are dealing with the partitioning of the whole
$\mathbb{R}^{3}$ into finite cells.

In what follows, we shall always assume that $\Pi:=(\pi
_{l})_{l\in \mathbb{N}}$ has a finite (partition) diameter
\begin{equation} \Delta\mathbf{x}=\Delta\mathbf{x}(\Pi)
  :=\sup_{l\in \mathbb{N}} D (\pi_{l}) <\infty. \label{dx}
\end{equation}

We are interested in the convergence properties of the solutions of
(\ref{discr3-BI}).

More specifically, we investigate the convergence properties of
\begin{equation} f_{\Delta t,\Delta
    \mathbf{x}}(t,\mathbf{x},\mathbf{v}
  ):=\sum_{j=1}^{[[\frac{T}{\Delta t}]]}\chi _{j}(t)
  \tilde{f}^{j}(\mathbf{x},\mathbf{v}), \quad 0 \leq t<T,
\label{approx-sol}
\end{equation} where $\chi_{j}$ is the indicator function of
the real interval $ [(j-1)  \Delta t,j  \Delta t)$.

In this paper, we suppose that (\ref{ec-bg}) concerns a gas with finite total
mass, hence physically meaningful solutions of the equation are elements of
the positive cone $\mathbb{L}_{+}^{1}$ of
$\mathbb{L}^{1}:=\mathbb{L}^{1}(\mathbb{R}^{3}\times
\mathbb{R}^{3};d\mathbf{x}d\mathbf{v})$ -- real. Therefore, our results will
be formulated in the sense of the convergence in $\mathbb{L}^{1}$.

Our basic hypotheses is that the collision kernel satisfies Grad's soft
cutoff condition \cite{Grad}. More precisely, throughout, we maintain the
following assumption:
\begin{assumption} There exist two constants $b_{0}>0$ and
  $0 \leq\lambda <2$ such that
\begin{equation}
  \int\limits_{\mathbb{S}^{2}}b(|\mathbf{v}-\mathbf{v}_{*}|,
  \boldsymbol{\omega})d\boldsymbol{\omega }\leq b_{0}
  {|\mathbf{v}-\mathbf{v}_{\ast}|}^{- \lambda}.\label{bo}
\end{equation}
\label{ass1}
\end{assumption}
Remark that the collision kernel for Maxwellian molecules
corresponds to the particular case $\lambda=0$ in
(\ref{bo}).

The main result of this paper shows that if $f(t) \in
\mathbb{L}_{+}^{1}$ is a solution to the Cauchy
problem~(\ref{ec-bg}), which decays in positions and
velocities at infinity, and has suitable Lipschitz
regularity properties (similar to those required in
\cite{BBb}), then $f_{\Delta t,\Delta \mathbf{x}}(t):=f_{\Delta
  t,\Delta \mathbf{x}}(t,\cdot, \cdot)
\in\mathbb{L}_{+}^{1}$ and $ f_{\Delta t,\Delta
  \mathbf{x}}(t) \overset{\mathbb{L}^{1}}\rightarrow f(t) $
as $(\Delta t+\Delta\mathbf{x})\rightarrow 0$, uniformly in
$\Delta t$ and $\Delta \mathbf{x}$. In addition, the
approximation $f_{\Delta t,\Delta \mathbf{x}}$ is consistent
with the laws of global conservation for mass, momentum and
energy, respectively,
\begin{equation} \int_{\mathbb{R}^{3}\times
    \mathbb{R}^{3}}d\mathbf{x}d\mathbf{v} \varphi
  _{i}(\mathbf{v})f_{\Delta t,\Delta
    \mathbf{x}}(t,\mathbf{x},\mathbf{v}
  )=\int_{\mathbb{R}^{3}\times
    \mathbb{R}^{3}}d\mathbf{x}d\mathbf{v} \varphi
  _{i}(\mathbf{v})f_{0}(\mathbf{x},\mathbf{v}),
\label{conservatiune}
\end{equation} where $\varphi _{i}(\mathbf{v)}$,
$i=0,1,2,3,4$ are as in (\ref{intJ}).

\section{Main result \label{s2}}
\subsection{Basic notations and definitions} \label{s3.1}
In general, in a given context of this paper, different
constants are differently denoted.  If some constant $c$
depends on parameters $\alpha_{1},\alpha_{2}, \dots$, we
will also denote it by $c_{\alpha_{1},\alpha_{2}, \dots}$, in
order to make explicit its parameter dependence.

Let $0\leq \alpha <\infty$, $0 <\tau < \infty $ and $m_{\alpha,
\tau}(\mathbf{x}, \mathbf{v}):=\exp(-\alpha {\mathbf{x}^{2}}-\tau
{\mathbf{v}^{2}})$.

We consider the following $\mathbb{L}^{\infty}$ - weighted
spaces. By $\mathbf{B}_{\tau }$, we denote the subspace of
elements $g\in \mathbb{L}^{\infty }:=\mathbb{L}^{\infty
}(\mathbb{R}^{3}\times \mathbb{R}^{3};d\mathbf{x}d
\mathbf{v})$ -- real satisfying $|| g|| _{\mathbf{B}_{\tau
  }}:=||m_{0,\tau}^{-1}g||_{\mathbb{L}^{\infty}} <\infty$.
Moreover, $\mathbf{M}_{\tau }$ denotes the subspace of the
elements $g\in\mathbf{B}_{\tau }$ with
$||g||_{\mathbf{M}_{\tau}}:=||m_{\tau,\tau}^{-1}g||_{\mathbb{L}^{\infty}}
<\infty$. Generically, we will refer to $||
\cdot||_{\mathbf{B}_{\tau}}$ and $|| \cdot
||_{\mathbf{M}_{\tau}}$ as the $\mathbf{B}$ - norm and
$\mathbf{M}$ - norm, respectively. We denote by
$\mathbf{B}_{\tau, + }$ and $\mathbf{M}_{\tau, + }$ the
positive cones of $\mathbf{B}_{\tau }$ and $\mathbf{M}_{\tau
}$, respectively, considered with the natural order $\leq$
(induced by $\mathbb{L}^{\infty}$).

Notice that (\ref{ut}) defines groups of linear isometries $\left\{
U^{t}\right\} _{t\in \mathbb{R}}$ on $\mathbf{B}_{\tau }$ and
$\mathbb{L}^{1}$,
\begin{equation}
\begin{array}{l} \displaystyle
|| U^{t}g|| _{\mathbf{B}_{\tau }}=|| g|| _{\mathbf{B}_{\tau }},
\quad \forall g\in \mathbf{B}_{\tau },\quad
  t\in \mathbb{R},
\\ \\  \displaystyle
 || U^{t}g|| _{\mathbb{L}^{1}}=|| g|| _{
    \mathbb{L}^{1}}, \quad \forall g\in \mathbb{L}^{1},\quad t\in
  \mathbb{R}.
 \end{array}
  \label{izo-X}
\end{equation}
One can easily check that  if $0<\tau_{1}, T < \infty$, then there is
some $\Theta=\Theta(T, \tau_{1})<\tau_{1}$, depending on $T$ and $\tau_{1}$,
such that for any $g \in \mathbf{M}_{\tau_{1}}$, we have
\begin{equation} || U^{t}g||_{\mathbf{M}_{\tau}} \leq
|| g||_{\mathbf{M}_{\tau_{1}}},
  \quad 0 \leq |t| \leq T, \quad 0\leq \tau <\Theta.
 \label{nonbd}
\end{equation}

One can also observe that  (\ref{pi-omo}) defines $\pi$ as a linear
contraction both in $\mathbf{B}_{\tau }$ and $\mathbb{L}^{1}$:
\begin{equation}
\begin{array}{l}
||\pi g||_{\mathbf{B}_{\tau }} \leq ||g||_{\mathbf{B}_{\tau }}, \quad \forall
g
\in
\mathbf{B}_{\tau },
\\  \\
||\pi g||_{\mathbb{L}^{1}} \leq ||g||_{\mathbb{L}^{1}}, \quad
\forall g \in \mathbb{L}^{1}.
\end{array}
\label{contr}
\end{equation}

Further, on $\mathbf{B}_{\tau }$, we consider the group of spatial
translations $\left\{T_{\mathbf{y}}\right\}_{\mathbf{y}\in
  \mathbb{R}^{3}}$ defined by
$(T_{\mathbf{y}}g)(\mathbf{x},\mathbf{v})
:=g(\mathbf{x}+\mathbf{y},\mathbf{v}) \label{ty}$ for almost
all $\mathbf{x}$, $\mathbf{v}$.  Obviously,
\begin{equation} || T_{\mathbf{y}}g|| _{\mathbf{B}_{\tau }}=||g|| _{
    \mathbf{B}_{\tau }}, \quad
\mathbf{y}\in\mathbb{R}^{3},
\label{izo-X1}
\end{equation}
and
\begin{equation} \label{uet}
(U^{t}g)(\mathbf{x},
\mathbf{v})=(T_{-(t\mathbf{v})}g)(\mathbf{x},\mathbf{v}),
\quad t\in \mathbb{R},
\end{equation}  for almost all $\mathbf{x}, \mathbf{v}\in
\mathbb{R}^{3}$.

Clearly, $U^{t}$, $\pi$ and $T_{\mathbf{y}}$ preserve the
order, in particular the positivity, in $\mathbf{B}_{\tau}$,
$\mathbf{M}_{\tau}$ and $\mathbb{L}^{1}$.

We need to introduce some suitable subsets of
$\mathbf{M}_{\tau}$ on which the map $\mathbf{y}\rightarrow
T_{\mathbf{y}}$ is Lipschitz continuous.

Let $0 < R,M < \infty $. By $\mathcal{M}_{\tau }(R,M)$, we
denote the subset of the elements $g\in \mathbf{M}_{\tau }$
satisfying
\begin{equation}
||g||_{\mathbf{M_{\tau}}}\leq R,
\label{ctau_R}
\end{equation} and
\begin{equation}
    ||(T_{\mathbf{y}}g)-g||_{\mathbf{M_{\tau}}} \leq M|\mathbf{y}|, \quad
 |\mathbf{y}|\leq 1.
\label{cor-liplip}
\end{equation}
It follows that $\mathcal{M}_{\tau }(R,M)\subset \mathbf{M}_{\tau } \subset
\mathbb{L}^{1}\cap \mathbf{B}_{\tau }$. Moreover, straightforward
calculations imply that, for any $ g \in \mathcal{M}_{\tau }(R,M)$,  we have
\begin{equation} ||g||_{\mathbf{B}_{\tau }}\leq
  ||g||_{\mathbf{M}_{\tau }}\leq R, \label{embed1}
\end{equation}
\begin{equation} ||g||_{ \mathbb{L}^{1}}\leq {\left (
      \frac{\pi}{\tau}\right )}^{3}
  ||g||_{\mathbf{M}_{\tau }}\leq {\left (
      \frac{\pi}{\tau}\right )}^{3} R,
\label{embed2}
\end{equation}
\begin{equation}
  ||\pi g - g||_{\mathbf{B}_{\tau}} \leq M \Delta \mathbf{x}, \quad
  \Delta \mathbf{x} \leq 1,
\label{dxinf}
\end{equation}
and
\begin{equation}
  ||\pi g - g||_{\mathbb{L}^{1}} \leq {\left
      (\frac{\pi}{\tau}\right )}^{3}M \Delta \mathbf{x},
  \quad \Delta \mathbf{x} \leq 1.
\label{dxlunu}
\end{equation}

We set $\mathcal{M}_{\tau, + }(R,M):=\{g\in\mathcal{M}_{\tau }(R,M): g \geq
0\}$.
\subsection{Main theorem and a practical criterion} \label{s3.2}
We are interested in approximating the $\mathbb{L}^{1}$ -
mild solutions of (\ref{ec-bg}). These are solutions of the
equation (see e.g., \cite{BePa})
\begin{equation}
  f(t)=U^{t}f_{0}+\int\limits_{0}^{t}U^{t-s}J(f(s))ds,
  \label{mild}
\end{equation}the integration with respect to $ds$ being in the sense of
Riemann in $\mathbb{L}^{1}$.  Here, the nonlinear operator $J$ is defined on
its natural domain of $\mathbb{L}^{1}$ by (\ref{jgh}) (with
$b(|\mathbf{v}-\mathbf{v}|, \boldsymbol{\omega})$ satisfying
assumption (\ref{bo})).

\begin{definition} Let $T>0$. An element $f\in
  C(0,T;\mathbb{L}^{1})$ is called mild solution on $[0,T]$
  for the Cauchy problem~(\ref{ec-bg}) in $\mathbb{L}^{1}$,
  if it satisfies Eq.~(\ref{mild}).
\label{defmild}
\end{definition}

Now we are in position to formulate our main result.

Let $0<\tau,T,R,M < \infty$, and denote by $\Lambda$ the couple of parameters
$(b_0,\lambda )$ appearing in Assumption \ref{ass1}.

Consider $f_{\Delta
t,\Delta \mathbf{x}}(t)$ as in (\ref{approx-sol}).

\begin{theorem} Suppose that the Cauchy problem~(\ref{ec-bg}) has a mild
solution $f(t)\in \mathcal{M}_{\tau, + }(R,M)$, $0 \leq t \leq T$. Then, for
each $0<\sigma <\tau $, there are some numbers $0<X_{*} < 1$, $0<T_{*} <
\min(1,T)$, and $ 0<K_{i} <\infty$, $i=1,2$, depending on $R,M,T,\tau,\sigma$
and $\Lambda $, such that, for any $0<\Delta t\leq T_{*}$ and $0<\Delta
\mathbf{x}\leq X_{*}$, one has
\begin{equation} f_{\Delta t,\Delta
\mathbf{x}}(t)\in \mathbb{L}_{+}^{1}\cap \mathbf{B}_{\sigma}, \quad 0\leq t <
T \label{pozitiv}
\end{equation}
and
\begin{equation} \sup_{t\in \lbrack 0,T]} || f_{\Delta
    t,\Delta \mathbf{x} }(t)-f(t)||_{{
      \mathbb{L}^{1}}}\leq K_{1}  \Delta t + K_{2} \Delta
  \mathbf{x} .
\label{convdiscr1}
\end{equation}
Moreover, $f_{\Delta t,\Delta \mathbf{x}}(t)$
satisfies property (\ref{conservatiune}) for all $0\leq t <
T$.
\label{conv-discr}
\end{theorem}

The proof of the theorem will be given in
Subsection~\ref{s4.2}. We only remark here that the proof
applies estimates based on the decomposition $J:=P-S$ into
the standard gain and loss operators \cite{Cer, BePa},
respectively, as well as estimates based on the similar
decomposition $J^{\pi}=P^{\pi}-S^{\pi}$, where
\begin{equation*}
P(g):=P_{B}(g,g), \quad S(g):=S_{B}(g,g)
\end{equation*}
and
\begin{equation*}
P^{\pi}(g):=P_{B}(g,\pi g), \quad S^{\pi}(g):=S_{B}(g,\pi g),
\end{equation*}
with $P_{B}$ and $S_{B}$ given by (\ref{ssb}). Due to assumption (\ref{bo}),
the above expressions define $P$, $S$, $P^{\pi}$ and $S^{\pi}$ as locally
Lipschitz continuous positive maps in $\mathbf{B}_{\tau }$ (in
$\mathbf{M}_{\tau }$).

A few remarks are  in order.

Since $P$, $S$ are locally Lipschitz continuous operators in
$\mathbf{B}_{\tau}$, and $\{U^{t}\}_{t \in \mathbb{R}}$ is a
continuous group of isometries on $\mathbf{B}_{\tau}$, the
contraction mapping principle applied to Eq.~(\ref{mild})
implies easily the existence and uniqueness of local-in-time
solutions (at least) in $\mathbf{B}_{\tau}$ (see, e.g.,
\cite{Grad} for the case of Maxwellian molecules). In fact,
there are more general results on the local and global
existence, uniqueness, and positivity of solutions to
Eq.~(\ref{mild}), in various spaces of functions decaying in
velocities and positions at infinity \cite{IlSh, BeTo, PO,
  BePa, GDc, GL, AlGa}.

We end this section with a sufficient condition for the
applicability of Theorem~\ref{conv-discr}, which allows for
replacing the Lipschitz assumption imposed in
Theorem~\ref{conv-discr} to the solutions of
Eq.~(\ref{mild}) by a similar one on the initial data.
\begin{theorem} Let $0 <\tau_{\ast}, R, T< \infty$. Suppose
  that the Cauchy problem~(\ref{ec-bg}) has a mild solution
  $f(t)\in \mathbf{M}_{\tau_{\ast},+}$ such that
\begin{equation}
||f(t)||_{\mathbf{M}_{\tau_{\ast}}} \leq R, \quad 0 \leq t \leq T.
\label{ftR}
\end{equation}
If there is a constant $M_{0}>0$ such that
\begin{equation}
f_{0} \in
\mathcal{M}_{\tau_{\ast}}(R,M_{0}),
\label{inico}
\end{equation}
then, for any $0<\tau_1<\tau_\ast$, there are two numbers,
$0<\tau<\tau_{\ast}$ (depending on $T,\tau_\ast$ and $\tau_1$) and $0<M<
\infty$ (depending on $R,M_0,T,\tau_\ast,\tau_1$ and $\Lambda $), such that
\begin{equation} f(t) \in \mathcal{M}_{\tau}(R,M), \quad
0\leq t \leq T.
\label{fcmrtau}
\end{equation}
\label{spatial-lip}
\end{theorem}
The proof of the theorem is given in Subsection \ref{s4.3}.

It should be observed that the subset of elements of
$\mathbf{M_{\tau_{\ast}}}$ satisfying (\ref{inico}) for some
$M_{0}>0$ is rather large, in the sense that
$\cup_{M_{0}>0}\mathcal{M}_{\tau_{\ast}}(R,M_{0})$ it is
dense in $\{g\in \mathbb{L}^{1}:
||g||_{\mathbf{M_{\tau_{\ast}}}}\leq R\}$ (with respect to
the topology of $\mathbb{L}^{1}$).
\section{Technical proofs \label{s3}}
\subsection{Auxiliary estimates} \label{s4.1}
In the following, some useful $\mathbf{B}$ - norm inequalities of \cite{BBb},
related to the properties of (\ref{discr3-BI}) and (\ref{mild}) are extended
to our setting. These are then supplemented with new $\mathbb{L}^{1}$ - norm
analogous inequalities. The main purpose is to prove Lemma~\ref{key} which
will play a key role in the next subsection, in estimating the cumulative
effect of the ``errors'' introduced when the solutions of the Boltzmann
equation are approximated by the solutions of (\ref{discr3-BI}).

The next lemma and its immediate consequences are based on standard
estimations on the collision operator, similar to those in \cite{BBb,BeTo,
AlGa,
  AlGa2}.  To formulate the lemma, first denote
\begin{equation}
G(\mathbf{v};\tau,\lambda):= \int\limits_{\mathbb{R}^{3}}
    |\mathbf{v}- \mathbf{v}_{\ast
  }|^{-\lambda} \exp(-\tau{{\mathbf{v}_{\ast}}}^{2}) d\mathbf{v}_{\ast};
  \quad 0\leq \lambda <2, \quad 0<\tau <\infty,
  \label{cgama2}
\end{equation}
and recall that by integrating upon $\mathbf{v}_{\ast}$, in
cylindrical coordinates in a frame with the $v_{\ast,3}$
axis along the direction of $\mathbf{v}$, one gets (see,
e.g., \cite{BeTo})
\begin{equation}
 sup_{\mathbf{v}\in
\mathbb{R}^{3}}G(\mathbf{v};\tau, \lambda) \leq
\pi^{\frac{3}{2}}\Pi\left (\frac{-\lambda}{2}\right
)\tau^{-\frac{3-\lambda}{2}},
\label{supgamma}
\end{equation}
where $\Pi(z)$ is the Gauss' Pi function.

With $b_{0}>0$  as in (\ref{bo}), set
\begin{equation}
  \begin{array} {l} \displaystyle
    b_{\Lambda}:=\Pi\left (\frac{-\lambda}{2}\right )b_{0},
    \\ \\ \displaystyle
    c_{\Lambda,\tau}=2\pi^{\frac{3}{2}}
    \Pi\left (\frac{-\lambda}{2}\right )\tau^{-\frac{3-\lambda}{2}}b_{0},
\end{array}
\label{cgama}
\end{equation}
as constants to be frequently used in the sequel.
\begin{lemma} a) For any
$m^{-1}_{\alpha,\tau}g_{i}, m^{-1}_{\alpha,\tau}h_{i} \in
\mathbb{L}^{\infty}$, $i=1,2$, one has
\begin{equation}
\begin{array}{l}\displaystyle
  ||m^{-1}_{\alpha,\tau}[J_{B}(g_{1},g_{2})-J_{B}(h_{1},h_{2})]||
  _{\mathbb{L}^{\infty}}
  \\ \\ \displaystyle
  \leq c_{\Lambda,\tau}( ||m^{-1}_{\alpha,\tau}
  g_{1}||_{\mathbb{L}^{\infty}}
  ||m^{-1}_{0,\tau}(g_{2}-h_{2})||_{\mathbb{L}^{\infty}}+
  ||m^{-1}_{\alpha,\tau} h_{2}||_{\mathbb{L}^{\infty}}
  ||m^{-1}_{0,\tau}(g_{1}-h_{1})||_{\mathbb{L}^{\infty}} ).
\end{array}
\label{jbinfty}
\end{equation}

b) For any $g_{i},h_{i} \in \mathbb{L}^{1}\cap \mathbf{B}_{\tau}$, $i=1,2$,
one has
\begin{equation}
||J_{B}(g_{1},g_{2})-J_{B}(h_{1},h_{2})||_{\mathbb{L}^{1}} \leq
c_{\Lambda,\tau}(||g_{1}||_{\mathbf{B}_{\tau}}
||g_{2}-h_{2}||_{\mathbb{L}^{1}}+ ||h_{2}||_{\mathbf{B}_{\tau}}
||g_{1}-h_{1}||_{\mathbb{L}^{1}}).
\label{jblunu}
\end{equation}
 \label{standard}
\end{lemma}
\textbf{Proof:} Since $J_{B}=P_{B}-S_{B}$, the proof follows
from suitable norm estimations of $S_{B}$ and $P_{B}$, respectively. To
this end, first observe that by applying  (\ref{bo}) in  (\ref{ssb}), one can
write
\begin{equation}
\begin{array}{l}
  |S_{B}(g_{1},g_{2})(\mathbf{x},\mathbf{v})-S_{B}
  (h_{1},h_{2})(\mathbf{x},\mathbf{v})|
  \\ \\ \displaystyle
  \leq b_{0}\int_{\mathbb{R}^{3}} {|\mathbf{v}-\mathbf{v}_{\ast
    }|}^{-\lambda}
  (|g_{1}|| g_{2\ast}-h_{2\ast}|+| h_{2\ast}||g_{1}-h_{1}|)
  d\mathbf{v}_{\ast}
  \end{array}
 \label {sg1g2}
\end{equation}
and
\begin{equation}
\begin{array}{l}
  |P_{B}(g_{1},g_{2})(\mathbf{x},\mathbf{v})
  -P_{B}(h_{1},h_{2})(\mathbf{x},\mathbf{v})|
  \\ \\ \displaystyle
  \leq b_{0}\int\limits_{\mathbb{R}^{3}}{|\mathbf{v}-\mathbf{v}_{\ast
    }|}^{-\lambda}
  (|g^{\prime }_{1}|| g^{\prime }_{2\ast}-h^{\prime }_{2\ast}|+| h^{\prime
  }_{2\ast}||g^{\prime }_{1}-h^{\prime }_{1}|)
  d\mathbf{v}_{\ast}.
\label {pg1g2}
\end{array}
\end{equation}
a) For the moment, let us denote the right hand side of (\ref{jbinfty}) by
\begin{equation*}
  {Q^{ h_{1},h_{2}}_{g_{1},g_{2}}}= c_{\Lambda,\tau}( ||m^{-1}_{\alpha,\tau}
  g_{1}||_{\mathbb{L}^{\infty}}
  ||m^{-1}_{0,\tau}(g_{2}-h_{2})||_{\mathbb{L}^{\infty}}+
  ||m^{-1}_{\alpha,\tau} h_{2}||_{\mathbb{L}^{\infty}}
  ||m^{-1}_{0,\tau}(g_{1}-h_{1})||_{\mathbb{L}^{\infty}} ).
\end{equation*}
In the integrand of (\ref{sg1g2}), we introduce  the obvious inequalities
\begin{equation*}
|g_{1}(\mathbf{x},\mathbf{v})|\leq {m}_{\alpha,
\tau}(\mathbf{x},\mathbf{v})|| {m}_{\alpha,
\tau}^{-1}g_{1}||_\mathbf{\mathbb{L^{\infty}}}, \quad
|h_{2}(\mathbf{x},\mathbf{v_{\ast}})|\leq {m}_{\alpha,
\tau}(\mathbf{x},\mathbf{v_{\ast}})||{m}_{\alpha,
\tau}^{-1}h_{2}||_\mathbf{\mathbb{L^{\infty}}},
\end{equation*}
\begin{equation*}
|g_{2}(\mathbf{x},\mathbf{v_{\ast}})-h_{2}(\mathbf{x},
\mathbf{v_{\ast}})|\leq {m}_{0,
\tau}(\mathbf{x},\mathbf{v_{\ast}})||{m}_{0,
\tau}^{-1}(g_{2}-h_{2})||_\mathbf{\mathbb{L^{\infty}}},
\end{equation*}
and
\begin{equation*}|g_{1}(\mathbf{x},\mathbf{v})-h_{1}(\mathbf{x},\mathbf{v})|\leq
{m}_{0, \tau}(\mathbf{x},\mathbf{v})||{m}_{0,
\tau}^{-1}(g_{1}-h_{1})||_\mathbf{\mathbb{L^{\infty}}}.
\end{equation*}
Then,  due to  (\ref{supgamma}), we find easily
\begin{equation*}
||m^{-1}_{\alpha,\tau}[S_{B}(g_{1},g_{2})-S_{B}(h_{1},h_{2})]||_{\mathbb{L}
^{\infty}} \leq \frac{{Q^{ h_{1},h_{2}}_{g_{1},g_{2}}}}{2}.
\end{equation*}
We proceed similarly with (\ref{pg1g2}), by taking advantage of the property
$\mathbf{v}^{\prime 2}+ \mathbf{v}^{\prime 2}_{\ast}= \mathbf{v}^{ 2}+
\mathbf{v}^{2}_{\ast}$ asserted in (\ref{econs}). We obtain
\begin{equation*}
||m^{-1}_{\alpha,\tau}[P_{B}(g_{1},g_{2})-P_{B}(h_{1},h_{2})]||_{\mathbb{L}
^{\infty}} \leq \frac{{Q^{ h_{1},h_{2}}_{g_{1},g_{2}}}}{2}.
\end{equation*}
This completes the proof of a).

b) We introduce inequalities
$|g_{1}(\mathbf{x},\mathbf{v})|\leq {m}_{0,
  \tau}(\mathbf{x},\mathbf{v})||g_{1}||_{\mathbf{B}_{\tau}}$
and $|h_{2}(\mathbf{x},\mathbf{v})|\leq {m}_{0,
  \tau}(\mathbf{x},\mathbf{v})||h_{2}||_{\mathbf{B}_{\tau}}$
in (\ref{sg1g2}).  Applying (\ref{supgamma}) to
estimate the integral over $\mathbf{v_{\ast}}$, and taking
the $\mathbb{L}^{1}$ - norm, we get
\begin{equation*}
  ||S_{B}(g_{1},g_{2})-S_{B}(h_{1},h_{2})||_{\mathbb{L}^{1}}
  \leq \frac {c_{\Lambda, \tau}}{2}( ||g_{1}||_{\mathbf{B}_{\tau}}
  ||g_{2}-h_{2}||_{\mathbb{L}^{1}}+ ||h_{2}||_{\mathbf{B}_{\tau}}
  ||g_{1}-h_{1}||_{\mathbb{L}^{1}} ).
\end{equation*}
Then, after the standard change of variables
$(\mathbf{v},\mathbf{v}_{\ast}) \rightarrow
(\mathbf{v}^{\prime}, \mathbf{v}^{\prime}_{\ast })$ in the
integral of (\ref{pg1g2}), similar computations as before
give
\begin{equation*}
  ||P_{B}(g_{1},g_{2})-P_{B}(h_{1},h_{2})||_{\mathbb{L}^{1}},
  \leq \frac{c_{\Lambda,\tau}}{2},
\end{equation*}
concluding  the proof of b). $\square$
\begin{corollary}
a) For any $g,h\in
  \mathbf{B}_{\tau }$, one has
\begin{equation}
|| J(g)-J(h)||_{\mathbf{B}_{\tau }}\leq
 c_{\Lambda,\tau}(|| g|| _{\mathbf{B}_{\tau }}+
 || h|| _{\mathbf{B}_{\tau }})  || g-h|| _{\mathbf{B}_{\tau }},
\label{JOd}
\end{equation}
\begin{equation} || J^{\pi }(g)-J^{\pi }(h)||
  _{\mathbf{B}_{\tau }}\leq c_{\Lambda,\tau}(|| g||
  _{\mathbf{B}_{\tau }}+|| h|| _{\mathbf{B}_{\tau
    }})  || g-h|| _{\mathbf{ B}_{\tau }}.
\label{J1-pi}
\end{equation}
Moreover, for any $g,h\in \mathbf{M}_{\tau}$,
\begin{equation} || J(g)-J(h)|| _{\mathbf{M}_{\tau }}\leq
  c_{\Lambda,\tau}(|| g|| _{\mathbf{M}_{\tau }}+|| h||
  _{\mathbf{M}_{\tau }})  || g-h||
  _{\mathbf{B}_{\tau }}.
\label{J00}
\end{equation}
b) Let $ 0 < \tau < \infty$. Then, for any $ g,h\in
  {\mathbb{L}^{1}}\cap \mathbf{B}_{\tau }$, one has
  \begin{equation} || J(g)-J(h)||_{{\mathbb{L}^{1}}}\leq
    c_{\Lambda,\tau}(|| g|| _{\mathbf{B}_{\tau }}+|| h||
    _{\mathbf{B} _{\tau }}) || g-h||
    _{{\mathbb{L}^{1}}}, \label{J2}
\end{equation}
\begin{equation} || J^{\pi }(g)-J^{\pi }(h)||
  _{{\mathbb{L}^{1}}}\leq c_{\Lambda, \tau}  (|| g||
  _{\mathbf{B}_{\tau }}+|| h|| _{ \mathbf{B}_{\tau
    }})  || g-h|| _{{\mathbb{L}^{1}}}. \label{J2-pi}
\end{equation}
\label{J-lip}
\end{corollary}
\textbf{Proof:} a) Let $g_{1}=g_{2}=h_{1}=h_{2}=g$ in (\ref{jbinfty}).  Then,
to obtain (\ref{JOd}) and (\ref{J00}), it is sufficient to  set $\alpha=0$
and $\alpha=\tau$, respectively, in (\ref{jbinfty}). Due to (\ref{jip0}),
inequality (\ref{J1-pi}) follows also from (\ref{jbinfty}) by setting
$g_{1}=g$, $g_{2}=\pi(g)$, $h_{1}=h$, $h_{2}=\pi(h)$, $\alpha=0$, and
applying property (\ref{contr}).

b) To obtain (\ref{J2}), we put $g_{1}=g_{2}=h_{1}=h_{2}=g$ in
(\ref{jblunu}).  Finally, due to (\ref{jip0}), we get (\ref{J2-pi}) by
setting $g_{1}=g$, $g_{2}=\pi(g)$, $h_{1}=h$, $h_{2}=\pi(h)$  in
(\ref{jbinfty}) and applying property (\ref{contr}). $\square $
\begin{remark} Since $J(0)=J^{\pi }(0)=0$, it follows that setting $h=0$ in
  each of the inequalities (\ref{JOd})--(\ref{J2-pi}) gives
  \begin{equation} ||J(g)||_{\mathbf{B}_{\tau }} \leq
    c_{\Lambda,\tau}||g||_{\mathbf{B}_{\tau }}^{2}, \quad
    ||J^{\pi}(g)||_{\mathbf{B}_{\tau }} \leq
    c_{\Lambda,\tau}||g||_{\mathbf{B}_{\tau }}^{2},
\label {gbtau}
\end{equation}
\begin{equation} ||J(g)||_{\mathbf{M}_{\tau }} \leq
  c_{\Lambda,\tau}||g||_{\mathbf{M}_{\tau }}
  ||g||_{\mathbf{B}_{\tau }},
 \label{gctau}
\end{equation}
\begin{equation} ||J(g)||_{\mathbb{L}^{1}} \leq
  c_{\Lambda,\tau}||g||_{\mathbf{B}_{\tau}}
  ||g||_{\mathbb{L}^{1}}, \quad
  ||J^{\pi}(g)||_{\mathbb{L}^{1}} \leq
  c_{\Lambda,\tau}||g||_{\mathbf{B}_{\tau
    }}||g||_{\mathbb{L}^{1}}.
 \label{gbtaul}
\end{equation}
\label{loclip}
\end{remark}

With $\pi $ and $\Delta x$ defined  in Section \ref{s1b}, we have:
\begin{lemma} For any $g\in \mathcal{M}_{\tau }(R,M)$, one has
\begin{equation} || J^{\pi}(g)-J(g)|| _{\mathbf{B}_{\tau
    }}\leq c_{\Lambda,\tau}MR
  \Delta\mathbf{x}, \quad \Delta\mathbf{x} \leq 1,
\label{J-pi}
\end{equation}
\begin{equation} || J^{\pi }(g)-J(g))||
  _{\mathbb{L}^{1}}\leq 2
  \pi^{\frac{9}{2}}b_{\Lambda}\tau^{-\frac{9-\lambda}{2}}MR
  \Delta\mathbf{x}, \quad \Delta\mathbf{x} \leq 1. \label{J-pidd}
\end{equation}
\label{J-lip-lem}
\end{lemma}
\textbf{Proof:} Due to (\ref{jip0}), we can set
$g_{1}=h_{1}=h_{2}=g$, $g_{2}=\pi g$ and $\alpha=0$ in
inequality (\ref{jbinfty}). This gives
\begin{equation*}
||J^{\pi}(g)-J(g)||_{\mathbf{B}_{\tau}} \leq c_{\Lambda,\tau}
||g||_{\mathbf{B}_{\tau}}||\pi g-g||_{\mathbf{B}_{\tau}}.
\end{equation*}
Similarly, by taking $g_{1}=h_{1}=h_{2}=g$, $g_{2}=\pi g$ in
inequality (\ref{jblunu}), we get
\begin{equation*} || J^{\pi }(g)-J(g)|| _{\mathbb{L}^{1}}
  \leq c_{\Lambda, \tau}  || g|| _{\mathbf{B}_{\tau }} || \pi g- g||_{
    \mathbb{L}^{1}}.
\end{equation*}
To conclude the proof,  recall that  $g$  satisfies (\ref{ctau_R}),
(\ref{dxinf}) and (\ref{dxlunu}). $\square $

Let $\max (x,y)$ denote the maximum between the real numbers $x$ and $y$.

\begin{lemma} For any $g\in \mathcal{M}_{\tau
}(R,M)$, one has
\begin{equation} ||U^{s}g-g|| _{\mathbf{B}_{\sigma }}\leq
  2^{\frac{1}{2}}e^{-\frac{1}{2}}{(\tau-\sigma)}^{-\frac{1}{2}}\max(R,M)|s|,
  \quad 0\leq \sigma <\tau,
\label{g-sigmadt}
\end{equation}
and
\begin{equation} || U^{s}g-g|| _{\mathbb{L}^{1}}\leq
4\pi^{\frac{5}{2}}\tau^{-\frac{7}{2}} \max(R,M)|s|.
\label{g-sigma}
\end{equation}
 \label{g-s-estim}
\end{lemma}
\textbf{Proof:} If $ |s\mathbf v| \leq 1$,  by applying (\ref{uet}) in
(\ref{cor-liplip}), we can write
\begin{equation*}
|(U^{s}g)(\mathbf{x},\mathbf{v})-
g(\mathbf{x},\mathbf{v})|\leq
M  |s|   |\mathbf v|   \exp(-\tau{\mathbf{v}}^{2})
\exp(-\tau \mathbf{x}^2).
\end{equation*}
Otherwise, we use (\ref{uet}) and ({\ref{ctau_R}}) in the
trivial inequality $|U^{s}g)(\mathbf{x},\mathbf{v})-
g(\mathbf{x},\mathbf{v})|\leq
|(U^{s}g)(\mathbf{x},\mathbf{v})|+
|g(\mathbf{x},\mathbf{v})|$, and  multiplying by $
|s\mathbf v| $, we get
\begin{equation*}
\begin{array}{l}
|(U^{s}g)(\mathbf{x},\mathbf{v})-
g(\mathbf{x},\mathbf{v})|
\\ \\
\leq R  |s|  |\mathbf v|   \exp(-\tau{\mathbf{v}}^{2})
[\exp[-\tau{(\mathbf{x}-s\mathbf{v})}^{2}+\exp(-\tau \mathbf{x}^2)], \quad
|s\mathbf v| > 1.
\end{array}
\end{equation*}
We can combine the above inequalities into
\begin{equation}
\begin{array}{l}
|(U^{s}g)(\mathbf{x},\mathbf{v})-
g(\mathbf{x},\mathbf{v})|
\\ \\       \displaystyle
\leq  \max(R,M)  |s|  |\mathbf v|   \exp(-\tau{\mathbf{v}}^{2})
[\exp[-\tau{(\mathbf{x}-s\mathbf{v})}^{2}+\exp(-\tau \mathbf{x}^2)],
\quad s\in\mathbb{R},
\end{array}
\label{usbl}
\end{equation}
for almost all $(\mathbf{x},\mathbf{v})$.  To obtain
(\ref{g-sigmadt}),  one takes the
$||\cdot||_{\mathbf{B}_{\sigma}}$ - norm of
(\ref{usbl}) and  applies
the inequality
$|\mathbf{v}|\exp(\sigma\mathbf{v}^2)\exp(-\tau\mathbf{v}^2)
\leq {[2e(\tau-\sigma)]}^{-\frac{1}{2}} $. Inequality (\ref{g-sigma}) follows
from the $\mathbb{L}^{1}$ - norm integration of
(\ref{usbl}). $\square$.

We can now prove the key result of this subsection.

Define
\begin{equation} I(g,h,t,s):=J^{\pi }(U^{t}g)-U^{s}J(h), \quad \forall
g,h\in\mathbf{B}_{\tau }, \quad t,s \in \mathbb{R}. \label{i}
\end{equation}
\begin{lemma} Let $0 <\sigma <\tau < \infty$. There are some constants
$k_{i}>0$, $i=1,2$ (depending on $R,M,\tau,\sigma$ and $\Lambda $), such that:

a) For any $g\in \mathbf{B}_{\sigma}$, $h,\hat{h} \in \mathcal{M}_{\tau
}(R,M)$ and $\Delta \mathbf{x} \leq 1$, one has
\begin{equation} ||I(g,h,t,s)||_{\mathbf{B}_{\sigma }}\leq
  c_{\Lambda, \sigma} (||g||_{\mathbf{B}_{\sigma }}+R)
  ||g-\hat{h} ||_{\mathbf{B}_{\sigma }}+k_{1} (||\hat{h}
  -h||_{\mathbf{B}_{\sigma }}+\Delta
  \mathbf{x}+|t|+|s|). \label{ij1}
\end{equation}

 b) For any $g\in \mathbb{L}^{1}\cap\mathbf{B}_{\sigma}$,
$h,\hat{h} \in \mathcal{M}_{\tau}(R,M)$ and $\Delta \mathbf{x} \leq 1$, one
has
\begin{equation} ||I(g,h,t,s)||_{\mathbb{L}^{1}}\leq
  c_{\Lambda, \sigma} (||g||_{\mathbf{B}_{\sigma }}+R)
  ||g-\hat{h} ||_{\mathbb{L}^{1}}+k_{2} (||\hat{h} -h||_{
    \mathbb{L}^1}+\Delta \mathbf{x}+|t|+|s|). \label{ij2}
\end{equation}
\label{key}
\end{lemma} \textbf{Proof:} We write
$I(g,h,t,s)=\sum\limits_{i=1}^{5}I_{i}$, where $
I_{1}=J^{\pi }(U^{t}g)-J^{\pi }(U^{t}\hat{h} )$,
$I_{2}=J^{\pi }(U^{t}\hat{h} )-J^{\pi }(\hat{h} )$,
$I_{3}=J^{\pi }(\hat{h} )-J(\hat{h} )$, $I_{4}=J(\hat{h}
)-U^{s}J(\hat{h} )$, and $I_{5}:=U^{s}J(\hat{h}
)-U^{s}J(h)$.

We proceed to establishing norm estimates for each $I_{i}$, keeping in mind
that, as elements of $\mathcal{M}_{\tau }(R,M)$, $h$ and $\hat{h}$ satisfy
(\ref{embed1}), and also belong to $\mathbf{B}_{\sigma }\cap \mathbb{L}^{1}$.

First, we apply inequality (\ref{J1-pi}) to $I_{1}$ and
$I_{2}$.  Then, making use of (\ref{izo-X}), we get for
$g\in \mathbf{B} _{\sigma }$,
\begin{equation} ||I_{1}||_{\mathbf{B}_{\sigma }}\leq
  c_{\Lambda, \sigma} (||g||_{\mathbf{B}_{\sigma }}+R)
  ||g-\hat{h} ||_{\mathbf{B}_{\sigma }}
 \label{i1b}
\end{equation} and
\begin{equation} ||I_{2}||_{\mathbf{B}_{\sigma }}\leq
  2c_{\Lambda, \sigma} R ||U^{t}\hat{h} -\hat{h}
  ||_{\mathbf{B}_{\sigma }}.
\label{i2b1}
\end{equation}
Moreover, starting from
(\ref{J2-pi}) applied to $I_{1}$ and $I_{2}$, and  using  again (\ref{izo-X}),
we obtain for $g\in \mathbf{B}_{\sigma }\cap \mathbb{L}^{1}$,
\begin{equation} ||I_{1}||_{\mathbb{L}^{1}}\leq c_{\Lambda, \sigma}
(||g||_{\mathbf{B}_{\sigma }}+R)  ||g-\hat{h} ||_{\mathbb{L}^{1}}
\label{i1l}
\end{equation} and
\begin{equation} ||I_{2}||_{\mathbb{L}^{1}}\leq 2c_{\Lambda, \sigma}
R  ||U^{t}\hat{h} -\hat{h} ||_{\mathbb{L}^{1}}.
\label{i2l1}
\end{equation}
Then, by applying (\ref{g-sigmadt}) to
(\ref{i2b1}), and (\ref{g-sigma}) to (\ref{i2l1}), respectively, we obtain
\begin{equation}
||I_{2}||_{\mathbf{B}_{\sigma }}\leq k_{1,2}|t|
\label{i2b}
\end{equation}
and
\begin{equation} ||I_{2}||_{\mathbb{L}^{1}}\leq k_{2,2} |t|, \label{i2l}
\end{equation}
respectively,
where
$k_{1,2}=2^{\frac{5}{2}}{\pi}^{\frac{3}{2}}e^{-\frac{1}{2}}
b_{\Lambda}{(\tau-\sigma}) ^{-\frac{1}{2}}\sigma^{-\frac{3-\lambda}{2}}
R\max(R,M)$ and $k_{2,2}={(2\pi)}^{4}b_{\Lambda}{\tau^{-\frac{7}{2}}}
{\sigma}^{-\frac{3-\lambda}{2}} R\max(R,M)$.

We estimate the norms of $I_{3}$, by means of  (\ref{J-pi}) and
(\ref{J-pidd}). We get
\begin{equation} ||I_{3}||_{\mathbf{B}_{\sigma }}
\leq k_{1,3}
  \Delta\mathbf{x}, \quad \Delta \mathbf{x} \leq 1,
\label{i3b}
\end{equation}
with $k_{1,3}=c_{\Lambda, \sigma}MR$, and
\begin{equation} ||I_{3}||_{\mathbb{L}^{1}}\leq k_{2,3}
  \Delta\mathbf{x}, \quad \Delta \mathbf{x} \leq 1,
\label{i3l}
\end{equation}
with $k_{2,3}=2\pi^{\frac{9}{2}}
b_{\Lambda}\tau^{-\frac{9-\lambda}{2}}MR$.

To estimate $I_{4}$, observe that since $\hat{h} \in
\mathcal{M}_{\tau}(R,M)$, then (\ref{gctau}) implies $||J(\hat{h}
)||_{\mathbf{M}_{\tau}}\leq c_{\Lambda, \tau} R^{2}$. Moreover, by
(\ref{J00}) and the commutation property $T_{\mathbf{y}}J(\hat{h}
)=J(T_{\mathbf{y}}\hat{h} )$, it follows that
\begin{equation*}
||T_{\mathbf{y}}J(\hat{h} )-J(\hat{h} )||_ {\mathbf{M}_{\tau }}
\leq c_{\Lambda, \tau}
(|| (T_{\mathbf{y}}\hat{h} || _{\mathbf{M}_{\tau }}+||\hat{h} ||
  _{\mathbf{M}_{\tau }})  ||T_{\mathbf{y}}\hat{h} -\hat{h} ||
  _{\mathbf{B}_{\tau }}.
\end{equation*}
Suppose that $|\mathbf{y}| \leq 1$. Then
$||T_{\mathbf{y}}\hat{h} -\hat{h} ||_{\mathbf{B}_{\tau }}
\leq M$, because of (\ref{cor-liplip}). But $||\hat{h}
||_{\mathbf{M}_{\tau }} \leq R$. Then, clearly,
$||T_{\mathbf{y}}\hat{h} ||_{\mathbf{M}_{\tau }}\leq
(R+M)$. It follows that
\begin{equation*}
  ||T_{\mathbf{y}}J(\hat{h} )-J(\hat{h} )||_ {\mathbf{M}_{\tau }}
  \leq c_{\Lambda, \tau}
  M  (2R+M)     |\mathbf{y}|, \quad
  |\mathbf{y}|\leq 1.
\end{equation*}
It appears that $J(\hat{h} )\in
\mathcal{M}_{\tau}(R_{\ast},M_{\ast})$ with
$R_{\ast}=c_{\Lambda, \tau} R^{2}$ and $M_{\ast}=c_{\Lambda,
  \tau} M (2R+M)$.  Therefore, we can apply Lemma~\ref
{g-s-estim} to obtain
\begin{equation} ||I_{4}||_{\mathbf{B}_{\sigma}}\leq k_{1,4}
|s|,
\label{i4b}
\end{equation}
with $k_{1,4}=2^{\frac{3}{2}}e^{-\frac{1}{2}}\pi^{\frac{3}{2}}
b_{\Lambda}\tau^{-\frac{3-\lambda}{2}} {(\tau-\sigma)}^{-\frac{1}{2}}
\max(R^{2},M (2R+M))$, and
\begin{equation} ||I_{4}||_{\mathbb{L}^{1}}\leq k_{2,4}
|s|.
\label{i4l}
\end{equation}
with $k_{2,4}=8\pi^{4}b_{\Lambda}\tau^
{-\frac{10-\lambda}{2}}\max(R^{2}, M(2R+M))$.

Finally, by means of (\ref{JOd}) and (\ref{izo-X}),
\begin{equation} ||I_{5}||_{\mathbf{B}_{\sigma }}\leq
k_{1,5}
||\hat{h} -h||_{ \mathbf{B}_{\sigma }},
 \label{i5b}
\end{equation}
while, by (\ref{J2}) and (\ref{izo-X}),
\begin{equation} ||I_{5}||_{\mathbb{L}^{1}}\leq k_{2,5}
  ||\hat{h} -h||_{\mathbb{L}^{1}}, \label{i5l}
\end{equation}
with $k_{1,5}=k_{2,5}=2c_{\Lambda, \sigma}  R$.

Inequalities (\ref{ij1}) and (\ref{ij2}) follow directly
from the above estimates, with $k_{i}= \max_{2\leq j\leq
  5}(k_{i,j})$, $i=1,2$. $ \square $

We end this section with a useful lemma on the mild solutions of
(\ref{ec-bg}).
\begin{lemma} Let $f(t)$ be a mild solution of
  Eq.~(\ref{ec-bg}), and $f(t)\in \mathcal{M} _{\tau }(R,M)$
  for all $0\leq t \leq T$.

a) Let $0<\sigma <\tau $. Then for each $t,s\in \lbrack 0,T]$,
\begin{equation} ||f(t)-f(s)||_{\mathbf{B}_{\sigma }}\leq
 d_{1}|t-s|,
 \label{t-lipd}
\end{equation}
where
$d_{1}=2^{\frac{1}{2}}e^{-\frac{1}{2}}{(\tau-\sigma)}^{-\frac{1}{2}}\max(R,M)
+c_{\Lambda, \sigma} R^{2}$.

b)For each $t,s\in \lbrack 0,T]$,
\begin{equation} ||f(t)-f(s)||_{\mathbb{L}^{1}}\leq
d_{2}|t-s|, \label{t-lip}
\end{equation}
where $d_{2}=4\pi^{\frac{5}{2}}\tau^{-\frac{7}{2}}\max(R,M)
+c_{\Lambda, \sigma}\left (\frac{\pi}{\sigma } \right )^3
R^{2}$.
 \label{T-lip}
\end{lemma}
\textbf{Proof:} We can suppose $t\geq s$. From  Eq.~(\ref{mild}),
\begin{equation} f(t)-f(s)=U^{t-s}f(s)-f(s)
  +\int\limits_s^tU^{t-u } J(f(u)) du.
\label{dif}
\end{equation}

a) Taking the $\mathbf{B} _{\sigma}$ norm in (\ref{dif}) and using
(\ref{izo-X}), we have
\begin{equation} ||f(t)-f(s)||_{\mathbf{B}_{\sigma}}\leq
  ||U^{t-s}f(s)-f(s)||_{\mathbf{B} _{\sigma
    }}+\int_{s}^{t}||J(f(u))||_{\mathbf{B}_{\sigma }}du.
\label{dif-D-tau}
\end{equation}
Further,  we apply (\ref{g-sigmadt}) to estimate the first term in
(\ref{dif-D-tau}). We obtain
\begin{equation} ||U^{t-s}f_{0}-f_{0}||_{\mathbf{B}_{\sigma
}}\leq  2^{\frac{1}{2}}e^{-\frac{1}{2}}{(\tau-\sigma)}^{-\frac{1}{2}}\max(R,M)
 (t-s).  \label{t1}
\end{equation}
Furthermore, we estimate the integral term of (\ref{dif-D-tau}) by means of
(\ref{gbtau}) and applying property (\ref{embed1}) to $f$. We obtain
\begin{equation}
\int\limits_{s}^{t}||J(f(u))||_{\mathbf{B}_{\sigma }}du \leq
c_{\Lambda, \sigma} R^{2} (t-s).  \label{t2}
\end{equation} Now (\ref{t-lipd}) results by
introducing (\ref{t1}) and (\ref{t2}) in (\ref{dif-D-tau}).

b) Taking the $\mathbb{L}^1$ norm in (\ref {dif}) and using (\ref{izo-X}), we
get
\begin{equation} ||f(t)-f(s)||_{\mathbb{L}^1}
\leq||U^{t-s}f(s)-f(s)||_{\mathbb{L}^1 }+
\int_s^t||J(f(u))||_{\mathbb{L}^1}du. \label{dif-L1}
\end{equation} Using (\ref{g-sigma}) to estimate the first term in
the r.h.s. of (\ref{dif-L1}), we get
\begin{equation} ||U^{t-s}f(s)-f(s)||_{\mathbb{L}^1 }\leq
4\pi^{\frac{5}{2}}\tau^{-\frac{7}{2}} \max(R,M)(t-s). \label{tt1}
\end{equation}  To estimate the second term in the r.h.s. of (\ref{dif-L1}),
we apply (\ref{gbtaul}) and the fact that $f$ satisfies both (\ref{embed1})
and (\ref{embed2}). We obtain
\begin{equation} \int_s^t||J(f(u))||_{\mathbb{L}^1 }du \leq
  c_{\Lambda, \sigma} \left (\frac{\pi}{\sigma }\right )^3
  R^{2} (t-s). \label{tt2}
\end{equation}

Finally, (\ref{t-lip}) is a consequence of (\ref{tt1}) and (\ref{tt2})
introduced in (\ref{dif-L1}). $\square $
\subsection{Uniform $\mathbf{B}$ - norm boundedness of
  $(\tilde{f}^{j})_{j}$. Proof of  Theorem~\ref{conv-discr}}
\label{s4.2} In the following,
we suppose that the assumptions of Theorem~\ref {conv-discr} are satisfied.
Thus, under the conditions of Theorem~\ref {conv-discr}, the mild solution
$f(t)$ of the Cauchy problem (\ref{ec-bg}) satisfies
\begin{equation} f(t)\in \mathcal{M}_{\tau }(R,M)\subset
\mathbf{B}_{\tau }\subset \mathbf{B}_{\sigma}
\label{fbctau}
\end{equation} and
\begin{equation} ||f(t)||_{\mathbf{B}_{\sigma }}\leq
||f(t)||_{\mathbf{B}_{\tau }}\leq ||f(t)||_{\mathbf{M}_{\tau
}}\leq R,\ \label{fbtau}
\end{equation} for all $ 0\leq t\leq T$, $0<\sigma <\tau $.

Put simply, the central (convergence) part of the argument behind
Theorem~\ref {conv-discr} consists in obtaining an appropriate uniform
$\mathbb{L}^{1}$ - norm estimate for the difference $\tilde{f}^{j}-f(t_{j})$,
where $\tilde{f} ^{j}$ is given by the recurrence (\ref{discr3-BI}), and
$f(t_{j})$ is the solution of Eq.~(\ref{mild}) at moment $t_{j}:=j \Delta t$,
$ j=0,1,\ldots,[[T/\Delta t]]$. To this end, one first needs to establish the
boundedness of the sequence $(\tilde{f}^{j})_{j}$ in a suitable $\mathbf{B}$
- norm. Technically,  the  proof of Theorem~\ref {conv-discr} applies
Lemma~\ref{key}. Thus, since $f_{0}$ satisfies (\ref{fbctau}), a
straightforward induction (based on the application of (\ref{izo-X}) and
(\ref{gbtau}) to the recurrence (\ref{discr3-BI})) implies \begin{equation}
\tilde{f}^{j}\in \mathbf{B}_{\sigma }. \label{fjbtau} \end{equation}
Consequently, (\ref{mild}) and (\ref{discr3-BI}) can be combined in a
well-defined (at least in $\mathbf{B}_{\sigma }$) expression \begin{equation}
\tilde{f}^{j}-f(t_{j})=U^{\Delta t}(\tilde{f}^{j-1}-f(t_{j-1}))+
\int_{t_{j-1}}^{t_{j}}[J^{\pi }(U^{\Delta t}\tilde{f}
^{j-1})-U^{t_{j}-u}J(f(u))]du \label{dif_e_ap} \end{equation} ($j=1,2,\ldots
,[[T/\Delta t]]$). In essence,  the central estimates of the proof of
Theorem~\ref {conv-discr} are obtained by applying  Lemma~\ref{key}  to
(\ref{dif_e_ap}), based on the immediate observation that the integrand of
(\ref{dif_e_ap}) satisfies \begin{equation} J^{\pi}(U^{\Delta t}
\tilde{f}^{j-1})-U^{t_{j}-u}J(f(u))=I(\tilde{f}^{j-1},f(u), \Delta
t,t_{j}-u), \label{JI} \end{equation} where $I$ is defined by  (\ref{i}).

In detail,  we first prove the following proposition which provides the
aforementioned $\mathbf{B}$ - norm boundedness of the sequence
$(\tilde{f}^{j})_{j}$ (and also yields the convergence of
 $(\tilde{f}^{j})_{j}$ with respect to the $\mathbf{B}$ - norm).

For some $0<\sigma < \tau$, define
\begin{equation} \tilde{F}^{j}:=\max
  \{||\tilde{f}^{0}||_{\mathbf{B}_{\sigma }},||\tilde{f}
  ^{1}||_{\mathbf{B}_{\sigma }},\ldots
  ,||\tilde{f}^{j-1}||_{\mathbf{B}_{\sigma}}\}, \quad
  j=1,2,\ldots,[[T/\Delta t]]. \label{fr}
\end{equation}
Due to (\ref{fjbtau}), obviously,  $\tilde{F}^{j}< \infty$,
$j=1,2,\ldots,[[T/\Delta t]]$.
\begin{proposition} a) There exists a strictly increasing
  continuous function $C(\cdot ):[0, \infty) \rightarrow (0, \infty)$
 (parameterized by $R,M,T,\tau ,\sigma$ and $\Lambda $) such
  that for any recurrence of the form (\ref{discr3-BI}),
  with time-step $0<\Delta t < T$ and cell-partition
  diameter $0<\Delta \mathbf{x}<1$,
\begin{equation}
||\tilde{f}^{j}-f(t_{j})||_{\mathbf{B}_{\sigma }}\leq
C(\tilde{F}^{j})  (\Delta t+\Delta \mathbf{x}), \quad
j=1,2,\ldots, [[T /\Delta t]].
 \label{ppp1-BI}
\end{equation}
b) There exist some numbers $\rho >0$, $0< X_{0} < 1$ and $ 0< T_{0} <\min
(1,T)$ (depending on $R,M,T,\tau,\sigma$ and $\Lambda$) such that, for any
recurrence (\ref{discr3-BI}) with $0<\Delta \mathbf{x} \leq X_{0}$ and $
0<\Delta t \leq T_{0}$,
\begin{equation}
||\tilde{f}^{j}-f(t_{j})||_{\mathbf{B}_{\sigma }}\leq
C(R+\rho)  (\Delta t+\Delta \mathbf{x})
\label{discr-marg-a-BI}
\end{equation} and
\begin{equation} ||\tilde{f}^{j}||_{\mathbf{B}_{\sigma
}}\leq R+\rho,
\label{discr-marg-a-BIbis}
\end{equation} $j=1,2,\ldots, [[T /\Delta t]]$.
\label{est-del-tx-BI}
\end{proposition}
\textbf{Proof:} a) As $\tilde{f}^{j}\in \mathbf{B}_{\sigma
}$, one can apply Lemma~\ref{key}~a), with $\hat{h}
=f(t_{j-1})$, to (\ref{JI}). The resulting inequality
contains the expression
$||f(t_{j-1})-f(u)||_{\mathbf{B}_{\sigma }}$ which is then
estimated by Lemma~\ref{T-lip}~a). One finds
\begin{equation}
\begin{array}{l} ||J^{\pi }(U^{\Delta
t}\tilde{f}^{j-1})-U^{t_{j}-u}J(f(u))||_{\mathbf{B}_{\sigma }}
\leq c_{\Lambda,\sigma} (||\tilde{f}^{j-1}||_{\mathbf{B}_{\sigma
}}+R)  ||\tilde{f}
^{j-1}-f(t_{j-1})||_{\mathbf{B}_{\sigma }} \\ \\ +k_{1}
\lbrack \Delta \mathbf{x}+\Delta t+(1+d_{1})
(t_{j}-u)],
\end{array}
\label{key1}
\end{equation} where the constants  $k_{1}$ and  $d_{1}$ are
given by Lemmas~\ref{key}~a) and \ref{T-lip}~a), respectively. By applying
(\ref{key1}) to estimate the ${\mathbf{B}_{\sigma}}$ - norm of
(\ref{dif_e_ap}), we obtain that there is some number $0<k<\infty$ (depending
on $R,M,\tau,\sigma$ and $\Lambda$) such that
\begin{equation}
\begin{array}{l}
||\tilde{f}^{j}-f(t_{j})||_{\mathbf{B}_{\sigma }}\leq
\lbrack
1+\varphi(||\tilde{f}^{j-1}||_{\mathbf{B}_{\sigma}})
\Delta t]  ||\tilde{f}^{j-1}-f(t_{j-1})||_{\mathbf{
B}_{\sigma }} \\\\ +k\Delta t  (\Delta
t+\Delta \mathbf{x}), \quad j=1,2,\ldots,[[T/\Delta t]],
\end{array}
\label{tri-BI_dt2}
\end{equation} where $\varphi (x):=c_{\Lambda, \sigma} (x+R)$.

Fix some $j^{\ast }=1,2,\ldots, [[T/\Delta t]]$. Due to
(\ref{fr}) and the monotonicity of $\varphi $, in
(\ref{tri-BI_dt2}), we can apply the inequality $\varphi
(\tilde{f}^{j})\leq \varphi ( \tilde{F}^{j^{\ast }})$,
$j=0,1,\ldots,j^{\ast}-1$. We are thus led to the following
simple Gronwall - type discrete scheme
\begin{equation}
\begin{array}{l}
||\tilde{f}^{j}-f(t_{j})||_{\mathbf{B}_{\sigma }}\leq
\lbrack 1+\varphi ( \tilde{F}^{j^{\ast }})  \Delta
t]  ||\tilde{f}^{j-1}-f(t_{j-1})||_{ \mathbf{B}_{\sigma
}} \\\\ +k\Delta t  (\Delta t+\Delta
\mathbf{x}), \quad j=1,2,\ldots,j^{\ast}.
\end{array}
\label{tri-BI_dt3}
\end{equation} As $\tilde{f}^{0}=f(0)$, by iterating
(\ref{tri-BI_dt3}), we get
\begin{equation*} ||\tilde{f}^{j^{\ast }}-f(t_{j^{\ast
}})||_{\mathbf{B}_{\sigma }}\leq \lbrack 1+\varphi
(\tilde{F}^{j^{\ast }})  \Delta t]^{j^{\ast }}
\frac{k(\Delta t+\Delta \mathbf{x})}{\varphi
(\tilde{F}^{{ j^{\ast }}})}.
\end{equation*} However,  $ j^{\ast} \leq T / \Delta t$. Consequently,
\begin{equation*}
\begin{array}{l} \displaystyle ||\tilde{f}^{j^{\ast
}}-f(t_{j^{\ast }})||_{\mathbf{B}_{\sigma }}\leq \lbrack
1+\varphi (\tilde{F}^{j^{\ast }})  \Delta t]^{\frac{T}{
\Delta t}}  \frac{k(\Delta t+\Delta
\mathbf{x})}{\varphi ( \tilde{F}^{{j^{\ast }}})}
\\\\\displaystyle \leq \frac{k\exp [(\varphi
(\tilde{F} ^{j^{\ast }})T]}{\varphi (\tilde{F}^{{j^{\ast
}}})}  (\Delta t+\Delta \mathbf{x}).
\end{array}
\end{equation*}
Since $\varphi(\tilde{F}^{{ j^{\ast }}})\geq c_{\Lambda, \sigma} R$, it
follows that
\begin{equation*} ||\tilde{f}^{j^{\ast }}-f(t_{j^{\ast
}})||_{\mathbf{B}_{\sigma }}\leq C(\tilde{F} ^{j^{\ast
}})  (\Delta t+\Delta \mathbf{x}),
\end{equation*}
with
\begin{equation} C(x):=\frac{k }
{c_{\Lambda,\sigma}   R}\exp
[c_{\Lambda, \sigma}  T(x+R)]>0,
\label{fc-BI}
\end{equation}
which is strictly increasing in $x$ on $[0,\infty)$.  This
concludes the proof of a), because $j^{*} \leq [[T/\Delta
t]]$ is arbitrary.  \medskip

b) From (\ref{fc-BI}) it follows that there is $ \rho
>0$ such that $ 0<{\rho}/{C(R+ \rho)}<\min (1,T)$.  Let
$0<X_{0} < {\rho}/{C(R+ \rho)}$ and $T_{0}= { \rho}/{C(R+
  \rho)}-X_{0}$. Therefore, for any $\Delta \mathbf{x}\leq
X_{0}$ and $\Delta t\leq $ $T_{0}$, we have
\begin{equation} 0<C(R+ \rho)(\Delta t+\Delta
\mathbf{x})\leq C(R+ \rho )(T_{0}+X_{0}) \leq
\rho. \label{in}
\end{equation} As $||\tilde{f}^{0}||_{\mathbf{B}_{\sigma
}}=||f_{0}||_{\mathbf{B} _{\sigma }}\leq R+ \rho$, we get $||\tilde{f}^{1}-
f(t_{1})||_{\mathbf{B} _{\sigma }} \leq C(R+\rho)(\Delta t + \Delta
\mathbf{x})\leq \rho $, by virtue of (\ref{ppp1-BI}) and (\ref{in}). Since
$f(t_{1})$ satisfies (\ref{fbtau}), it follows that $||\tilde
f_{1}||_{\mathbf{B} _{\sigma }}\leq R+ \rho$. Then a straightforward
induction concludes the proof of b). $\square$

Based on Proposition~\ref{est-del-tx-BI}, we can now prove
the main result of the paper.

\textbf{Proof of Theorem~\ref{conv-discr}:} With the
notations of the theorem, let $0<\sigma <\tau$.

 \underline { a) Proof of (\ref{pozitiv}):} First we show that  $f_{\Delta
 t,\Delta
\mathbf{x}}(t)\in \mathbb{L}^{1}\cap \mathbf{B}_{\sigma }$.

By virtue of (\ref{fjbtau}), we need only prove that
\begin{equation} \tilde{f}^{j}\in \mathbb{L}^{1},\quad
j=0,1,\ldots ,[[T/\Delta t]].
\label{fjl1}
\end{equation}
Obviously, $\tilde{f}^{0}=f_{0}\in \mathcal{M}_{\tau
}(R,M)\subset \mathbb{L}^{1}$.  To check that
$\tilde{f}^{1}\in \mathbb{L}^{1}$, first observe that from
(\ref{dif_e_ap}), particularized to $j=1$, one finds
\begin{equation} ||\tilde{f}^{1}||_{\mathbb{L}^{1}}\leq
||f(t_{1})||_{\mathbb{L} ^{1}}+\int_{0}^{\Delta t}(||J^{\pi
}(U^{\Delta t}f_{0})||_{\mathbb{L} ^{1}}+||U^{\Delta
t-u}J(f(u))||_{\mathbb{L}^{1}})du. \label{dif_e_ap2}
\end{equation}
Since $f(t)\in\mathcal{M}_{\tau }(R,M)\subset
\mathbb{L}^{1}$, we need only check that the integral term
of the above inequality is finite.  To this
end, we apply (\ref{gbtaul}) and (\ref{izo-X}), to estimate
the terms of the sum under the integral in
(\ref{dif_e_ap2}). We get $||J^{\pi }(U^{\Delta
  t}f_{0})||_{\mathbb{L}^{1}}\leq c_{\Lambda, \sigma}
||f_{0}||_{ \mathbb{B}_{\sigma }}
||f_{0}||_{\mathbb{L}^{1}}$ and $||U^{\Delta
  t-u}J(f(u))||_{\mathbb{L}^{1}}\leq c_{\Lambda, \sigma}
||f(u)||_{\mathbf{B}_{\sigma }}
||f(u)||_{\mathbb{L}^{1}}$. Now it remains to observe that
$||f(u)||_{\mathbb{L}^{1}}$ satisfies (\ref{embed2}), and
that $||f(u)||_{\mathbf{B}_{\sigma }}\leq R$, by virtue of
(\ref{fbtau}).

As $\tilde{f}^{j}$ satisfies (\ref{discr-marg-a-BIbis}), the
proof of (\ref {fjl1}) is completed by induction, following
a similar argument as before, based on the application of
(\ref{gbtaul}) and (\ref{izo-X}).

To conclude the proof of (\ref{pozitiv}), it remains to show
that $\tilde{f}^j \geq 0$. We proceed by induction, applying
a trick as in \cite{Bab89}.

1) $f^0\geq 0$ by hypothesis.

2) By (\ref{discr3-BI}),
\begin{equation*} \tilde{f}^{j}=U^{\Delta
t}\tilde{f}^{j-1}-\Delta t  S^{\pi }(U^{\Delta t}
\tilde{f}^{j-1})+\Delta t  P^{\pi }(U^{\Delta
t}\tilde{f}^{j-1}), \quad j=1,\ldots,[[T/\Delta t]].
\end{equation*} Suppose that $\tilde{f}^{j-1}\geq 0$. As
$U^{\Delta t}$ and $P^{\pi }$ are positivity preserving
operators, in order to show that $f^{j}\geq 0$, it is
sufficient to prove that
\begin{equation} U^{\Delta t}\tilde{f}^{j-1}-\Delta t
S^{\pi }(U^{\Delta t}\tilde{f} ^{j-1})\geq 0. \label{positj}
\end{equation} Observe that
\begin{equation*} U^{\Delta t}\tilde{f}^{j-1}-\Delta t
S^{\pi }(U^{\Delta t}\tilde{f} ^{j-1})=U^{\Delta
t}\tilde{f}^{j-1}[1-\Delta t  E(U^{\Delta t}\tilde{f}
^{j-1})],
\end{equation*} where
\begin{equation*}
 E(g)(\mathbf{x},\mathbf{v}):=\sum_{l\in
\mathbb{N}}\chi _{l}( \mathbf{x})\frac{1}{| \pi
_{l}| }\int\limits_{\pi _{l}}d
\mathbf{y}\int_{\mathbb{R}^{3}\times
\mathbb{S}^{2}}d\mathbf{v}_{\ast }d\boldsymbol{\omega}
 b(|\mathbf{v}-\mathbf{v}_{\ast }|,\boldsymbol{\omega})
g(\mathbf{y},\mathbf{v}_{\ast}).
\end{equation*}
But (\ref{discr-marg-a-BIbis}) gives
\begin{equation*} (U^{\Delta t}\tilde{f}^{j-1})(\mathbf{y},\mathbf{v}_{\ast})
<(R+\rho)  \exp
(-\sigma |\mathbf{v_{\ast }} |^{2}), \quad
j=1,\ldots,[[T/\Delta t]],
\end{equation*}
for $0< \Delta \mathbf{x} \leq X_{0}$ and $0< \Delta t \leq
T_{0}$ as in Proposition~\ref{est-del-tx-BI}. Then, by
virtue of (\ref{supgamma}),
\begin{equation*} E(U^{\Delta t}\tilde{f}^{j-1})\leq
\frac{1}{2}c_{\Lambda, \sigma}(R+\rho)=:\mathcal{D}_{0},
\quad j=1,\ldots,[[T/\Delta t]].
\end{equation*} Therefore, it is sufficient to set $T_{*}=
\min (T_{0}, \mathcal{D}_{0}^{-1})$, in order that the inequality
(\ref{positj}) be satisfied.

\underline { b) Proof of (\ref{convdiscr1}):} We show that there is some
number $K>0$ (depending on  $R,M, T,\tau ,\sigma$ and $\Lambda $) such that
for all $j=0,1,\ldots,[[T/\Delta t]]$,
\begin{equation}
||\tilde{f}^{j}-f(t_{j})||_{\mathbb{L}^{1}}\leq K
(\Delta t+\Delta \mathbf{x}), \quad 0<\Delta t\leq
T_{*},\quad 0<\Delta \mathbf{x}\leq X_{*}, \label{prop2k-BI}
\end{equation} with $X_{*}= X_{0}$, where $X_{0}$ is as in
Proposition~\ref{est-del-tx-BI}.

We start as in the proof Proposition~\ref{est-del-tx-BI} a). Since $\tilde{f}
^{j}\in \mathbb{L}^{1}\cap \mathbf{B}_{\sigma }$, we apply
Lemma~\ref{key}~b), with $\hat{h} =f(t_{j-1})$, to (\ref{JI}).  The resulting
inequality contains the expression $||f(t_{j-1})-f(u)||_{\mathbb{L}^{1}}$
which is then estimated by Lemma~\ref{T-lip} b). We obtain
\begin{equation}
  \begin{array}{l} \displaystyle
    ||J^{\pi}(U^{\Delta t}
    \tilde{f}^{j-1})-U^{t_{j}-u}J(f(u))||_{\mathbb{L}^{1}} \leq
    c_{\Lambda,\sigma}(||\tilde{f}^{j-1}||_{\mathbf{B}_{\sigma }}+R)
    \\ \\ \displaystyle
    \times
    ||\tilde{f} ^{j-1}-f(t_{j-1})||_{\mathbb{L}^{1}}
    +k_{2}  \lbrack \Delta \mathbf{x} +\Delta t+(1+d_{2})  (t_{j}-u)],
\end{array}
\label{key2}
\end{equation} where the constants $k_{2}$ and $d_{2}$ are
given by Lemmas~\ref{key}~b) and \ref{T-lip}~b), respectively. We apply
(\ref{key2}) to estimate the $\mathbb{L}^{1}$ - norm of (\ref{dif_e_ap}).
Also, we take advantage of (\ref{izo-X}) and of the key property
(\ref{discr-marg-a-BIbis}). After a straightforward computation, it follows
that there exist two constants $\tilde{K}_{i}>0$, $i=1,2$ (which depend on
$R,M,T,\tau,\sigma$ and $\Lambda$)  such that
\begin{equation}
\begin{array}{l}
  ||\tilde{f}^{j}-f(t_{j})||_{\mathbb{L}^{1}}\leq
  (1+\tilde{K}_{1}  \Delta t)
  ||\tilde{f}^{j-1}-f(t_{j-1})||_{\mathbb{L}^{1}} \\\\
  +\tilde{K}_{2}  \Delta t  (\Delta t+\Delta \mathbf{x}),
  \quad j=1,2,\ldots,[[T/\Delta t]].
\end{array}
\label{convL-discr3b}
\end{equation} The simple scheme (\ref{convL-discr3b}) can
be iterated directly with respect to
$||\tilde{f}^{j}-f(t_{j})||_{\mathbb{L}^{1}}$. As
$\tilde{f}^{0}=f(0)$, we get
\begin{equation*}
  ||\tilde{f}^{j}-f(t_{j})||_{\mathbb{L}^{1}}\leq \frac{\tilde{K}
    _{2}}{\tilde{K}_{1}}  (1+\tilde{K}_{1}\Delta t)^{j}  (\Delta
  t+\Delta \mathbf{x}).
\end{equation*} However,  $j \leq [[T/\Delta t]]$. Therefore,
\begin{equation*}
  ||\tilde{f}^{j}-f(t_{j})||_{\mathbb{L}^{1}} \leq \frac{\tilde{K}_{2}}
  {\tilde{K}_{1}}  (1+\tilde{K}_{1}\Delta t)^{\frac{T}{\Delta
      t}}  (\Delta t+\Delta \mathbf{x})\leq \frac{\tilde{K}_{2}}{\tilde{K}_{1}}
        \exp (\tilde{K}_{1}T)  (\Delta t+\Delta \mathbf{x}),
\end{equation*} for all $j=1,2,\ldots,[[T/\Delta
t]]$. Finally, set
\begin{equation*} K=\frac{\tilde{K}_{2}}
{\tilde{K}_{1}}  \exp (\tilde{K}_{1}T),
\end{equation*} we obtain (\ref{prop2k-BI}).

Now (\ref{convdiscr1}) follows directly  from (\ref{prop2k-BI}) and
(\ref{t-lip}),  with $K_{1}=d_{2}+K$ and $K_{2}=K$, where $d_{2}$ is the
constant  of (\ref{t-lip}).

\underline { c)  Proof of (\ref{conservatiune}):} The conservation property
(\ref{conservatiune}) follows from (\ref{intJ}) and
(\ref{discr3-BI}), by observing that, by virtue of
(\ref{jgh}) and (\ref{pi-omo}), one obtains
\begin{equation*}
  \begin{array}{l} \displaystyle
\int_{\mathbb{R}^{3}\times\mathbb{R}^{3}}\varphi_{i}(\mathbf{v})J_{B}(U^{\Delta
t}\tilde{f}^{j-1},\pi U^{\Delta
t}\tilde{f}^{j-1})(\mathbf{x},\mathbf{v})d\mathbf{x}d\mathbf{v}
\\ \displaystyle
=\int_{\mathbb{R}^{3}\times\mathbb{R}^{3}}\varphi_{i}(\mathbf{v})\pi
J_{B}(U^{\Delta t}\tilde{f}^{j-1},\pi U^{\Delta
t}\tilde{f}^{j-1})(\mathbf{x},\mathbf{v})d\mathbf{x}d\mathbf{v}
\\ \displaystyle =
\int_{\mathbb{R}^{3}\times\mathbb{R}^{3}}\varphi_{i}(\mathbf{v})J_{B}(\pi
U^{\Delta t}\tilde{f}^{j-1},\pi U^{\Delta
t}\tilde{f}^{j-1})(\mathbf{x},\mathbf{v})d\mathbf{x}d\mathbf{v},
\quad i=0,1,2,3,4,
\end{array}
\end{equation*} with $\varphi_{i}(\mathbf{v})$ as in
(\ref{intJ}). Then it is sufficient to invoke (\ref{intJ}). $ \square $

\subsection{Proof of Theorem~\ref{spatial-lip}} \label{s4.3}
Since $||f(t)||_{\mathbf{M}_{\tau}} \leq
||f(t)||_{\mathbf{M}_{\tau{\ast}}} \leq R$ for all $0 \leq
t\leq T$, $0< \tau \leq \tau_{\ast}$, it is sufficient to
prove that $f(t)$ satisfies an inequality of the form
(\ref{cor-liplip}).

We start with the remark that since $T_{\mathbf{y}}$
commutes with $U^{t}$ and $J$, then, by virtue of
(\ref{mild}), we have, for any $0\leq t\leq T$,
\begin{equation}
T_{\mathbf{y}}f(t)-f(t)=U^{t}(T_{\mathbf{y}}f_{0}-f_{0})+\int
\limits_{0}^{t}U^{t-s}[J(T_{\mathbf{y}}f(s))-J(f(s))]ds.
\label{tyf}
\end{equation}
Let $0<\tau_{1}<\tau_{\ast}$. Due to (\ref{nonbd}), there is
some $0<\tau <\tau_{1}$, such that
\begin{equation}
||T_{\mathbf{y}}f(t)-f(t)||_{\mathbf{M}_{\tau}} \leq
||T_{\mathbf{y}}f_{0}-f_{0}||_{\mathbf{M}_{\tau_{1} }}
 +\int\limits_{0}^{t}||J(T_{\mathbf{y}}f(s))-J(f(s))]||
_{\mathbf{M}_{\tau_{1}}}ds.
\label{tyctauini}
\end{equation}
Observing that the first term of the sum in (\ref{tyctauini}) satisfies
\begin{equation*}
  ||(T_{\mathbf{y}}f_{0}-f_{0})||_{\mathbf{M}_{\tau_{1} }} \leq
  || (T_{\mathbf{y}}f_{0}-f_{0})||_{\mathbf{M}_{\tau_{*}}}\leq
  M_{0}  |\mathbf{y}|, \quad |\mathbf {y}|\leq 1,
\end{equation*}
and   introducing (\ref{J00}) in the integral term of (\ref{tyctauini}),  we
obtain
\begin{equation}
  \begin{array}{l} \displaystyle
    ||T_{\mathbf{y}}f(t)-f(t)||_{\mathbf{M}_{\tau}} \leq
    M_{0}  |\mathbf{y}|
    \\ \\ \displaystyle
    +c_{\Lambda,\tau_{1}}
    \int\limits_{0}^{t}(||T_{\mathbf{y}}f(s)||_{\mathbf{M}_{\tau_{1} }} +
    ||f(s)||_{\mathbf{M }_{\tau_{1}
      }})||T_{\mathbf{y}}f(s)-f(s)||_{\mathbf{B}_{\tau_{1}}}ds.
\end{array}
\label{tyctau}
\end{equation}
Further, we estimate the factors of the product under the
integral sign in (\ref{tyctau}). To this end, by observing
that a straightforward computation gives
$||T_{\mathbf{y}}f(s)||_{\mathbf{M}_{\tau_{1}}}\leq R
\exp\left ({{\frac
      {\tau_{\ast}^2}{\tau_{\ast}-\tau_{1}}}}x\right )$, for
all $0 \leq s \leq T $, and $|\mathbf{y}| \leq 1$, and using
$||f(s)||_{\mathbf{M}_{\tau_{1}}} \leq
||f(s)||_{\mathbf{M}_{\tau_{\ast}}} \leq R$, we get
\begin{equation}
||T_{\mathbf{y}}f(s)||_{\mathbf{M}_{\tau_{1} }} +
||f(s)||_{\mathbf{M }_{\tau_{1}
}} \leq R\left [1+\exp\left ({{\frac
{\tau_{\ast}^2}{\tau_{\ast}-\tau_{1}}}}\right )\right ], \quad |\mathbf{y}| \leq 1,
\quad 0\leq s\leq T.
\label{fterm}
\end{equation}
To estimate the second factor of the product under the
integral sign of (\ref{tyctau}), first observe that
\begin{equation} ||T_{\mathbf{y}}f_{0}-f_{0}|| _{\mathbf{B}_{\tau_{1}}} \leq
M_{0} |\mathbf{y}|, \quad |\mathbf{y}| \leq 1, \label{lip-0b}
\end{equation}
because of the assumption $f_{0} \in
\mathbf{M}_{\tau_{\ast}}(R,M_{0})$.  Then, due to property
(\ref{lip-0b}), a standard argument applied to (\ref{tyf})
gives
 \begin{equation} || T_{\mathbf{y}}f(t)-f(t)||
_{\mathbf{B}_{\tau_{1} }}\leq M_{0}  \exp (2c_{\Lambda,\tau_{1}}Rt)
|\mathbf{y}|,\quad |\mathbf{y}| \leq 1, \quad 0\leq t\leq T.
\label{tym0}
\end{equation}
Indeed, by taking the $\mathbf{B}_{\tau_{1}}$ - norm of
(\ref{tyf}), and applying (\ref{izo-X}), (\ref{lip-0b}) and
(\ref{JOd}), we obtain
\begin{equation*}
\begin{array}{l}\displaystyle
  ||T_{\mathbf{y}}f(t)-f(t)|| _{\mathbf{B}_{\tau_{1} }}\leq M_{0}
  |\mathbf{y}|
  \\ \\ \displaystyle
  +c_{\Lambda,\tau_{1}} \int\limits_{0}^{t}(||
  T_{\mathbf{y}}f(s)||_{\mathbf{B}_{\tau_{1} }}+|| f(s)||
  _{\mathbf{B}_{\tau_{1} }})  ||T_{\mathbf{y}}f(s)-f(s)||
  _{\mathbf{B}_{\tau_{1} }}ds, \quad |\mathbf{y}| \leq 1.
\end{array}
\end{equation*}
But  (\ref{izo-X1}) and (\ref{ftR}) imply $
||T_{\mathbf{y}}f(t)||_{\mathbf{B}_{\tau_{1}}}=||f(t)||_{\mathbf{B_{\tau_{1}}}}
\leq ||f(t)||_{\mathbf{M}_{\tau_{\ast}}} \leq R$, hence
\begin{equation*}
  ||T_{\mathbf{y}}f(t)-f(t)||_{\mathbf{B}_{\tau_{1}}}\leq M_{0}
  |\mathbf{y}|+2c_{\Lambda,\tau_{1}}  R  \int\limits_{0}^{t} ||
  T_{\mathbf{y}}f(s)-f(s)||_{\mathbf{B}_{\tau_{1} }}ds,
\end{equation*} so that the application of Gronwall's inequality yields
(\ref{tym0}).

Thus, by using (\ref{fterm}) and (\ref{tym0}) in
(\ref{tyctau}), we finally obtain\begin{equation*}
  ||T_{\mathbf{y}}f(t)-f(t)||_{\mathbf{M}_{\tau}}\leq M
  |\mathbf{y}|, \quad |\mathbf{y}| \leq 1, \quad 0\leq t
  \leq T,
\end{equation*}
with
\begin{equation}
  M=M_{0}
  \left \{1+\frac{1}{2}(\exp
    (2c_{\Lambda,\tau_{1}}RT)-1)\left [1+\exp\left ({{\frac
        {\tau_{\ast}^{2}}{\tau_{\ast}-\tau_{1}}}}\right )\right ]\right \},
\label{m}
\end{equation}
hence (\ref{fcmrtau}) is satisfied with $M$  of (\ref{m}). This concludes the
proof. $\square$
\section{Example and conclusions \label{s5}}
We present a simple application of Theorem~\ref{spatial-lip} to the solutions
of the Cauchy problem for the Boltzmann equation near vacuum. (We skip over
mentioning parametric dependencies as not being particularly relevant for our
purposes.)

In what follows, it is sufficient to consider the existence and uniqueness of
local in time, positive solutions to Eq.~(\ref{mild}) for a small initial
datum bounded by a space-velocity Maxwellian. The following result can be
easily obtained by applying the Kaniel-Shinbrot monotone iteration scheme
\cite{BeTo, KaSh}, or by elementary fixed point methods \cite{GDc}.

\begin{proposition} Let $0<\tau_{0} < \infty$ and $0 \leq
  f_{0}\in \mathbf{M}_{\tau _{0} }$. For each $T>0$, there
  are $0 <r <R< \infty$ such that if $|| f_{0}||
  _{\mathbf{M}_{\tau _{0 }}}\leq r$, then Eq.~(\ref{mild})
  has a unique solution satisfying
\begin{equation}
0 \leq f(t,\mathbf{x}, \mathbf{v})\leq R
\exp[-\tau_{0}{(\mathbf{x}-t\mathbf{v})}^2-\tau_{0}{\mathbf{v}}^{2}],\quad 0
\leq t \leq T,
 \label{U-t}
\end{equation}
for almost all $(\mathbf{x}, \mathbf{v})\in \mathbb{R}^3 \times
\mathbb{R}^3$. \label{existence}
\end{proposition}
 Our application follows by combining  the above proposition
 with Theorem~\ref{spatial-lip}.
\begin{proposition}
  Let $0<\tau_{0}, M_{0} < \infty$. For each $T>0$, there
  are $0<r, R , M <\infty $ and $0<\tau <\tau_{0}$ such that
  if $0 \leq f_{0} \in \mathcal{M}_{\tau _{0 }}(r,M_{0})$,
  then the Cauchy problem~(\ref{ec-bg}) has a unique mild
  solution $0 \leq f(t)\in \mathcal{M}_{\tau}(R,M)$ for all
  $0\leq t \leq T$.
\end{proposition}
\textbf{Proof:} Consider the solution $f$ of
Eq.~(\ref{mild}) provided by the above proposition. By
virtue of (\ref{nonbd}), there is
$0<\Theta=\Theta(T,\tau_{0})<\tau_{0}$ such that for any
$0\leq \tau_{\ast} <\Theta$, one can write $||f(t)||_{
  \mathbf{M}_{\tau_{\ast}}}=||U^{t}U^{-t}f(t)||_{\mathbf{M}_{\tau_{\ast}}}
\leq ||U^{-t}f(t)||_{\mathbf{M}_{\tau_0}}$, $0\leq t\leq
T$. Then, due to (\ref{U-t}),
\begin{equation} || f(t)||_{\mathbf{M}_{\tau_{\ast}}}\leq
R,\quad 0\leq t\leq T.
\label{f-sig}
\end{equation}
Obviously, $f(t)\in \mathbb{L}^{1}$, for all $t \geq
0$. Moreover, $f_{0} \in \mathcal{M}_{\tau _{\ast
  }}(R,M_{0})$, because $ f_{0} \in \mathcal{M}_{\tau _{0
  }}(r,M_{0})$, $0<\tau_{\ast}<\tau_{0}$, and $r<R$.
Consequently, Theorem~\ref{spatial-lip} applies, concluding
the proof.  $\square$

We end this sections with a few comments about our results.

Theorem~\ref{conv-discr} implies immediately the convergence
in discrepancy of (\ref{approx-sol}). Thus, the validation
of Nanbu's simulation scheme for the Boltzmann equation in
the whole space and Maxwellian molecular interactions can be
supplemented with a similar result as in \cite{BBb}.

To better clarify why the approach of \cite{BBb} is not directly applicable to
our setting (for the Boltzmann equation in the whole space), recall that
expression (\ref{discr-marg-a-BI}) establishes the convergence of the
solutions of the recurrence (\ref{discr3-BI}) in a suitable $\mathbf{B}$ -
norm, i.e. in a (Maxwellian-weighted) $\mathbb{L}^{\infty}$ - space. Such a
property was sufficient to ensure the convergence in $\mathbb{L}^{1}$, in the
setting of \cite{BBb} for the Boltzmann gas in a finite domain $\Omega$,
because of the continuous embedding of the space $\{h:\;m^{-1}_{0,\tau}h\in
\mathbb{L}^\infty(\Omega\times\mathbb{R}^3;d{\mathbf{x}}d\mathbf{v})\}$ into
$ \mathbb{L}^{1}(\Omega\times \mathbb{R}^{3};d{\mathbf{ x}}d{\mathbf{v}})$,
which holds when $\Omega$ is bounded. However, this is not the case if
$\Omega=\mathbb{R}^{3}$, when solely the weak convergence property
(\ref{discr-marg-a-BI}) is not  sufficient to imply the $\mathbb{L}^{1}$ -
convergence.

Theorem~\ref{spatial-lip} extends somehow the simpler property, mentioned in
  \cite{BBb}, that if the initial condition of Eq.~(\ref{mild}) satisfies
  $f_{0} \in \mathcal{B}_{\tau}(R,M)$, then the solution of the equation also
  satisfies $f(t) \in \mathcal{B}_{\tau}(R,M)$ for all $0 < t \leq T$.
  Nevertheless, the latter property  remains valid in the context of the
  Boltzmann equation in the entire space, being actually
  established within the proof of Theorem~\ref{spatial-lip}, by deriving
  inequality (\ref{tym0}) as a consequence of (\ref{lip-0b}).

  Following a line of reasoning as in the present paper,
  Theorem~\ref{conv-discr} can be generalized to a wider
  class of solutions of the Boltzmann equation, with slower
  decay at infinity, like those considered in some
  investigations on the Cauchy problem for the Boltzmann
  equation with near-vacuum conditions \cite{PO, BePa}.

  The results of this paper can be also extended to more
  complicated Boltzmann like models as those describing
  several spaces of chemically interacting fluids \cite{GDc,
    GDa, GDi}. A potential application would be the
  validation of the space-dependent Nanbu scheme for the
  reacting gas, by extending results obtained in the
  space-homogeneous case in \cite{GDd}, \cite{DM}.

  Due to the explicit form of the constants involved in the
  technical inequalities of Section~\ref{s3}, the proofs of
  Proposition~\ref{est-del-tx-BI} and
  Theorem~\ref{spatial-lip} may be also detailed to provide
  explicit upper bounds for the constants $K_1$ and $K_2$ of
  (\ref{convdiscr1}), as well as for the other constants
  appearing in Theorem~\ref{spatial-lip}. Such upper bounds
  may be useful in estimating the errors introduced by the
  approximation (\ref{discr3-BI}), as well as in the
  parametric optimization of the approximation. However, a
  detailed computation of the above bounds is beyond the
  scope of this paper.

\section*{Acknowledgment}

This work was partially supported by a grant of the Romanian
National Authority for Scientific Research, CNCS UEFISCDI,
project number PN-II-RU-TE-2012-3-0196

\end{document}